\documentclass{article}

\usepackage{arxiv}

\usepackage{amssymb}
\usepackage{amsmath}
\usepackage{amsthm}
\usepackage{bm}
\usepackage{natbib}
\usepackage{hyperref}
\usepackage{tikz}
\usepackage{tkz-euclide}
\usetkzobj{all} 
\usepackage{chngcntr}
\usepackage{verbatim}
\usepackage{mathtools}
\usepackage{esvect}
\usepackage{subfiles}
\usepackage{color}
\usepackage{xspace}
\usepackage{xparse}
\usepackage[noend]{algpseudocode}
\RequirePackage{algorithm}
\usepackage{graphicx}
\usepackage{subcaption}

\newcommand{\inner}[1]{\langle #1 \rangle}
\newcommand{\repeatn}[3]{\underbrace{\text{#1} \ldots \text{#2}}_{#3}}

\newtheorem{example}{Example}
\newtheorem{theorem}{Theorem}
\newtheorem{lemma}[theorem]{Lemma}
\newtheorem{proposition}[theorem]{Proposition}
\newtheorem{corollary}[theorem]{Corollary}
\newtheorem{remark}[theorem]{Remark}
\newtheorem{assumption}{Assumption}
\newtheorem{definition}{Definition}

\newcommand{\Oc}{\mathcal{O}}
\newcommand{\reals}{\mathbb{R}}
\newcommand{\gd}{\textsc{Gd}\xspace}
\newcommand{\nag}{\textsc{Nag}\xspace}
\NewDocumentCommand{\InternalVector}{mmm}{%
	\begin{bmatrix}#1\\ #2\\ #3\end{bmatrix}%
}
\NewDocumentCommand{\vektor}{>{\SplitArgument{2}{,}}m}{\InternalVector#1}

\author{\name{Jingzhao Zhang}\email{jzhzhang@mit.edu}\\
	\name{Aryan Mokhtari}\email{aryanm@mit.edu}\\
	\name{Suvrit Sra}\email{suvrit@mit.edu}\\
	\name{Ali Jadbabaie}\email{jadbabai@mit.edu}\\
	\addr{Laboratory for Information and Decision Systems\\
		Institute for Data, Systems, and Society\\
		Massachusetts Institute of Technology}
}

\renewcommand{\comment}[1]{}
\newcolumntype{S}{>{\centering\arraybackslash} m{.10\linewidth} }
\newcolumntype{T}{>{\centering\arraybackslash} m{.30\linewidth} }
\title{
Direct Runge-Kutta Discretization Achieves Acceleration
}
\begin{document}

\maketitle

\begin{abstract}
  We study gradient-based optimization methods obtained by directly discretizing a second-order ordinary differential equation (ODE) related to the continuous limit of Nesterov's accelerated gradient method. When the function is smooth enough, we show that acceleration can be achieved by a stable discretization of this ODE using standard Runge-Kutta integrators. Specifically, we prove that under Lipschitz-gradient, convexity and order-$(s+2)$ differentiability assumptions, the sequence of iterates generated by discretizing the proposed second-order ODE converges to the optimal solution at a rate of $\mathcal{O}({N^{-2\frac{s}{s+1}}})$, where $s$ is the order of the Runge-Kutta numerical integrator. Furthermore, we introduce a new local flatness condition on the objective, under which rates even faster than $\mathcal{O}(N^{-2})$ can be achieved with low-order integrators and only gradient information. Notably, this flatness condition is satisfied by several standard loss functions used in machine learning. We provide numerical experiments that verify the theoretical rates predicted by our results.


\end{abstract}

\section{Introduction}\label{sec:introduction}
In this paper, we study accelerated first-order optimization algorithms for the problem
\begin{equation}
  \label{main_problem}
  \min_{x \in \reals^d}\quad f(x),
\end{equation}
where $f$ is convex and sufficiently smooth. A classical method for solving~\eqref{main_problem} is gradient descent (\gd), which displays a sub-optimal convergence rate of $\Oc(N^{-1})$---i.e., the gap $f(x_N)-f(x^*)$ between \gd and the optimal value $f(x^*)$ decreases to zero at the rate of $\mathcal{O}(N^{-1})$. Nesterov's seminal accelerated gradient method~\citep{nesterov-smooth-acceleration} matches the oracle lower bound of $O(N^{-2})$~\citep{nemirovskii-lowerbound}, and is thus a central result in the theory of convex optimization.

However, ever since its introduction, acceleration has remained somewhat mysterious, especially because Nesterov's original derivation relies on elegant but unintuitive algebraic arguments. This lack of understanding has spurred a variety of recent attempts to uncover the rationale behind the phenomenon of acceleration~\citep{allen-linear,bubeck-geometric,lessard-iqc,hu-dissipativity,scieur-regularized,fazlyab-dynamical}.

We pursue instead an approach to \nag (and accelerated methods in general) via a continuous-time perspective. This view was recently studied by \citet{su-differential}, who showed that the continuous limit of \nag is a second order ODE describing a physical system with vanishing friction;  \citet{wibisono-variational} generalized this idea and proposed a class of ODEs by minimizing Bregman Lagrangians. 

Although these works succeed in providing a richer understanding of Nesterov's scheme via its continuous time ODE, they fail to provide a general discretization procedure that generates \emph{provably convergent} accelerated methods. In contrast, we introduce a second-order ODE that generates an accelerated first-order method for smooth functions if we simply discretize it using \emph{any} Runge-Kutta numerical integrator and choose a suitable step size. 

\subsection{Summary of results}
Assuming that the objective function is convex and sufficiently smooth, we establish the following:
\begin{list}{{\tiny$\blacksquare~$}}{\leftmargin=1.4em}
\vspace*{-4pt}
\item We propose a second-order ODE, and show that the sequence of iterates generated by discretizing using a Runge-Kutta integrator converges to the optimal solution at the rate $\Oc({N^{\frac{-2s}{s+1}}})$, where $s$ is the order of the  integrator. By using a more precise numerical integrator, (i.e., a larger $s$), this rate approaches the optimal rate $\Oc(N^{-2})$. 
\item We introduce a new local flatness condition for the objective function (Assumption~\ref{assump:growth}), under which Runge-Kutta discretization obtains convergence rates even faster than $\Oc(N^{-2})$, without requiring high-order integrators. In particular, we show that if the objective is  locally flat  around a minimum, by using \emph{only gradient information} we can obtain a convergence rate of $\Oc(N^{-p})$, where $p$ quantifies the degree of local flatness. Acceleration due to local flatness may seem counterintuitive at first, but our analysis reveals why it helps.
\end{list}
To the best of our knowledge, this work presents the first direct\footnote{That is, discretize the ODE with known numerical integration schemes without resorting to reverse engineering \nag's updates.} discretization of an ODE that yields accelerated gradient methods. Unlike~\citet{betancourt2018symplectic} who study symplecticity and consider variational integrators, and \citet{scieur-integration} who study consistency of integrators, we focus on the order of integrators (see \S\ref{sec:RK}). We argue that the stability inherent to the ODE and order conditions on the integrators suffice to achieve acceleration. 

\vspace*{-4pt}
\subsection{Additional related work}
\vspace*{-4pt}
Several works \citep{alvarez2000minimizing,attouch2000heavy,bruck1975asymptotic,attouch1996dynamical} have studied the asymptotic behavior of solutions to dissipative dynamical systems. However, these works retain a theoretical focus as they remain in the continuous time domain and do not discuss the key issue, namely,  stability of discretization. Other works such as \citep{krichene-accelerated}, study the counterpart of \citet{su-differential}'s work for mirror descent algorithms and achieve acceleration via Nesterov's technique. \citet{diakonikolas2017approximate} proposes a framework to analyze the first order mirror descent algorithms by studying ODEs derived from duality gaps. Also, \citet{raginsky2012} obtain nonasymptotic rates for continuous time mirror descent in a stochastic setting. 

A textbook treatment of numerical integration is given in \citep{hairer-textbook}; some of our proofs build on material from Chapters 3 and 9. \citep{isaacson1994analysis} and \citep{west2004variational} also provide nice introductions to numerical analysis.


\section{Problem setup and background}
Throughout the paper we assume that the objective $f$ is convex and sufficiently smooth. Our key result rests on two key assumptions introduced below. The first assumption is a \emph{local} flatness condition on $f$ around a minimum; our second assumption requires $f$ to have bounded higher order derivatives. These assumptions are sufficient to achieve acceleration simply by discretizing suitable ODEs without either resorting to  reverse engineering to obtain discretizations or resorting to other more involved integration mechanisms. 

We will require our assumptions to hold on a suitable subset of $\reals^d$. Let $x_0$ be the initial point to our proposed iterative algorithm. First consider the sublevel set
\begin{equation}\label{eq:set-S}
\mathcal{S} :=\{x \in \reals^d \mid f(x) \le \exp(1)((f(x_0) - f(x^*) + \|x_0 - x^*\|^2) + 1 \},
\end{equation}
where $x^*$ is a minimum of~\eqref{main_problem}. Later we will show that the sequence of iterates obtained from discretizing a suitable ODE never escapes this sublevel set. Thus, the assumptions that we introduce need to hold only within a subset of $\reals^d$. Let this subset be defined as
\begin{equation}\label{eq:set-A}
	\mathcal{A} :=\{x \in \reals^d \mid \exists x' \in \mathcal{S},\ \|x-x'\| \le 1\}, 
\end{equation}
that is, the set of points at unit distance to the initial sublevel set~\eqref{eq:set-S}. The choice of unit distance is arbitrary, and one can scale that to any desired constant.

\begin{assumption} \label{assump:growth} 
  There exists an integer $p \ge 2$ and a positive constant $L$ such that for any point $x\in \mathcal{A}$, and for all indices $i \in \{1,...,p-1\}$, we have the lower-bound
  \begin{equation}
    \label{eq:1}
    f(x) - f(x^*) \ge \tfrac{1}{L} \|\nabla^{(i)} f(x)\|^{\frac{p}{p-i}},
  \end{equation}
  where $x^*$ minimizes $f$ and $\|\nabla^{(i)}f(x)\|$ denotes the operator norm of the tensor $\nabla^{(i)}f(x)$.
\end{assumption}
Assumption~\ref{assump:growth} bounds high order derivatives by function suboptimality, so that these derivatives vanish as the suboptimality converges to $0$. Thus, it quantifies the flatness of the objective around a minimum.\footnote{One could view this as an error bound condition that reverses the gradient-based upper bounds on suboptimality stipulated by the Polyak-\L{}ojasiewicz condition~\citep{lojasiewicz-ensembles,attouch-proximal}.} When $p=2$, Assumption~\ref{assump:growth} is slightly weaker than the usual Lipschitz-continuity of gradients (see Example \ref{ex:p2}) typically assumed in the analysis of first-order methods, including \nag. If we further know that the objective’s Taylor expansion around an optimum does not have low order terms, p would be the degree of the first nonzero term.

\begin{example} \label{ex:p2}
  Let $f$ be convex with $\frac{L}{2}$-Lipschitz continuous gradients, i.e., $\|\nabla f(x)-\nabla f(y)\|\le \frac{L}{2}\|x-y\|$. Then, for any $x, y\in \reals^d$ we have
  \begin{equation*}
    f(x) \ge f(y) + \inner{\nabla f(y), x - y} + \tfrac{1}{L}\|\nabla f(x) - \nabla f(y)\|^2.
  \end{equation*}
  In particular, for $y=x^*$, an optimum point, we have $\nabla f(y)=0$, and thus we have $f(x)-f(x^*) \ge \tfrac{1}{L} \|\nabla f(x)\|^2,$ which is nothing but inequality~\eqref{eq:1} for $p=2$ and $i=1$.
\end{example}

\begin{example} \label{ex:4th}
  Consider the $\ell_p$-norm regression problem: $\min_x f(x) =  \|A x - b\|^p_p$, for even integer $p \ge 2$. If $\exists x^*, Ax^*=b$, then $f$ satisfies inequality \eqref{eq:1} for $p$, and $L$ depends on $p$ and the operator norm of $A$. 
\end{example}
Logistic loss satisfies a slightly different version of Assumption~\ref{assump:growth} because its minimum can be at infinity. We will explain this point in more detail in Section~\ref{sec:logistic}. 

Next, we introduce our second assumption that adds additional restrictions on differentiability and bounds the growth of derivatives.
\begin{assumption} \label{assump:differentiable}
  There exists an integer $s \ge p$ and a constant $M \ge 0$, such that $f(x)$ is order $(s+2)$ differentiable. Furthermore, for any $x\in \mathcal{A}$, the following operator norm bounds hold:
  \begin{equation}
    \label{eq:2}
    \|\nabla^{(i)} f(x)\| \le M, \qquad \text{for}\ i=p,p+1,\ldots,s,s+1,s+2.
  \end{equation}
\end{assumption}
When the sublevel sets of $f$ are compact and hence the set $\mathcal{A}$ is also compact; as a result, the bound~\eqref{eq:2} on high order derivatives is implied by continuity. In addition, an $L_p$ loss of the form $\|Ax-b\|_p^p$ also satisfy~\eqref{eq:2} with $M = p! \|A\|_2^p$.
 

\subsection{Runge-Kutta integrators} \label{sec:RK}
Before moving onto our new results (\S\ref{sec:main}) let us briefly recall \emph{explicit} Runge-Kutta (RK) integrators used in our work. For a more in depth discussion please see the textbook~\citep{hairer-textbook}.
\begin{definition} \label{def:rk}
  Given a dynamical system $\dot{y} = F(y)$, let the current point be $y_0$ and the step size be $h$. An explicit $S$ stage Runge-Kutta method generates the next step via the following update:
  \begin{align}\label{eq:runge-kutta}    
    g_i &= y_0 + h \sum_{j=1}^{i-1} a_{ij} F(g_j), \qquad
    \Phi_h(y_0) = y_0 + h \sum_{i=1}^{S} b_i F(g_i),
  \end{align}
where $a_{ij}$ and $b_i$ are suitable coefficients defined by the integrator;  $\Phi_h(y_0)$ is the estimation of the state after time step $h$, while $g_i$ (for $i=1,\ldots,S$) are a few neighboring points where the gradient information $F(g_i)$ is evaluated. 
\end{definition}
By combining the gradients at several evaluation points, the integrator can achieve higher precision by matching up Taylor expansion coefficients. Let $\varphi_h(y_0)$ be the true solution to the ODE with initial condition $y_0$; we say that an integrator $\Phi_h(y_0)$ has order $s$ if its \emph{discretization error} shrinks as
\begin{equation}
\label{eq:disc-error}
\|\Phi_h(y_0) - \varphi_h(y_0)\| = O(h^{s+1}), \qquad\text{as}\ h\to 0.
\end{equation}
In general, RK methods offer a powerful class of numerical integrators, encompassing several basic schemes. The \emph{explicit Euler's} method defined by $\Phi_h(y_0) = y_0 + hF(y_0)$ is an explicit RK method of order 1, while the \emph{midpoint} method $\Phi_h(y_0) = y_0 + hF(y_0+\tfrac{h}{2}F(y_0))$ is of order 2. Some high-order RK methods are summarized in \citep{highorderRK}. An order 4 RK method requires 4 stages, i.e., 4 gradient evaluations, while an order 9 method requires 16 stages.  


\section{Main results}\label{sec:main}
In this section, we introduce a second-order ODE and use explicit RK integrators to generate iterates that converge to the optimal solution at a rate faster than $\mathcal{O}(1/t)$ (where $t$ denotes the time variable in the ODE). A central outcome of our result is that, at least for objective functions that are smooth enough, it is not the integrator type that is the key ingredient of acceleration, but a careful analysis of the dynamics with a more powerful Lyapunov function that achieves the desired result. More specifically, we will show that by carefully exploiting boundedness of higher order derivatives, we can achieve both stability and acceleration at the same time.  

We start with Nesterov's accelerated gradient (\nag) method that is defined according to the updates
\begin{equation}
\label{eq:nesterov}
\begin{split}
x_k &= y_{k-1} - h\nabla f(y_{k-1}), \qquad y_k = x_k + \tfrac{k-1}{k+2}(x_k - x_{k-1}).
\end{split}
\end{equation}
\citet{su-differential} showed that the iteration~\eqref{eq:nesterov} in the limit is equivalent to the following ODE
\begin{align}\label{eq:nag-ode}
\ddot{x}(t) + \tfrac{3}{t}\dot{x}(t) + \nabla f(x(t)) = 0,  \qquad \text{where}\ \dot{x} = \tfrac{d x}{d t}
\end{align}
when one drives the step size $h$ to zero.  It can be further shown that in the continuous domain the function value $f(x(t))$ decreases at the rate of $\Oc({1}/{t^2})$ along the trajectories of the ODE. This convergence rate can be accelerated to an arbitrary rate in continuous time via time dilation as in ~\citep{wibisono-variational}. In particular, the solution to
\begin{align}\label{eq:wib-ode}
\ddot{x}(t) + \tfrac{p+1}{t}\dot{x}(t) + p^2 t^{p-2} \nabla f(x(t)) = 0,
\end{align}
has a convergence rate $\Oc({1}/{t^p})$. When $p>2$, \citet{wibisono-variational} proposed rate matching algorithms via utilizing higher order derivatives (e.g., Hessians). In this work, we focus purely on first-order methods and study the stability of discretizing the ODE directly when $p\ge2$. 

Though deriving the ODE from the algorithm is solved, deriving the update of \nag or any other accelerated method by directly discretizing an ODE is not. As stated in \citep{wibisono-variational}, explicit Euler discretization of the ODE in \eqref{eq:nag-ode} may not lead to a stable algorithm. Recently, \citet{betancourt2018symplectic} observed empirically that Verlet integration is stable and suggested that the stability relates to the symplectic property of the Verlet integration. However, in our proof, we found that \textbf{\emph{the order condition of Verlet integration would suffice to achieve acceleration}}. Though symplectic integrators are known to preserve modified Hamiltonians for dynamical systems, we weren't able to leverage this property for dissipative systems such as~\eqref{eq:ode}.

This principal point of departure from previous works underlies Algorithm 1, which solves~\eqref{main_problem} by discretizing the following ODE with an order-$s$ integrator:
\begin{align}\label{eq:ode}
\ddot{x}(t) + \frac{2p+1}{t}\dot{x}(t) + p^2t^{p-2} \nabla f(x(t)) = 0.
\end{align}
where we have augmented the state with time, to turn the non-autonomous dynamical system into an autonomous one. The solution to \eqref{eq:ode} exists and is unique when $t>0$. This claim follows by local Lipschitzness of $f$ and is discussed in more details in Appendix A.2 of \cite{wibisono-variational}.
\begin{algorithm}[t] \label{algorithm_main}
	Algorithm 1: Input($f, x_0, p, L, M, s, N$)\Comment{Constants $p, L, M$ are the same as in Assumptions}
	\begin{algorithmic}[1]
		\State Set the initial state $y_0 = [\vec{0}; x_0; 1]\in \reals^{2d+1}$
		\State Set step size h = $C/N^{\frac{1}{s+1}}$ \Comment{C is determined by $p, L, M, s, x_0$}
		\State $x_N \gets \text{Order-s-Runge-Kutta-Integrator}(F, y_0, N, h)$ \Comment{F is defined in equation \ref{eq:dynamics}}
		\State \textbf{return} $x_N$
	\end{algorithmic}
\end{algorithm}

We further highlight that the ODE in \eqref{eq:ode} can also be written as the dynamical system
\begin{align}\label{eq:dynamics}
\dot{y} = 
F(y)
= 
\begin{bmatrix}
-\frac{2p+1}{t} v - p^2t^{p-2} \nabla f(x)\\
v\\
1
\end{bmatrix},\qquad\text{where}\ y = [v; x; t].
\end{align}
We have augmented the state with time to obtain an autonomous system, which can be readily solved numerically with a Runge-Kutta integrator as in Algorithm 1. To avoid singularity at $t=0$, Algorithm 1 discretizes the ODE starting from $t=1$ with initial condition $y(1) = y_0 =  [0; x_0; 1]$. The choice of $1$ can be replaced by any arbitrary positive constant. 

Notice that the ODE in \eqref{eq:ode} is slightly different from the one in \eqref{eq:wib-ode}; it has a coefficient $\frac{2p+1}{t}$ for $\dot{x}(t)$ instead of $\frac{p+1}{t}$. This modification is crucial for our analysis via Lyapunov functions (more details in Section \ref{sec:proof} and Appendix~\ref{supp_prop_dec_energy}).  

The parameter $p$ in the ODE~\eqref{eq:ode} is set to be the same as the constant in Assumption \ref{assump:growth} to achieve the best theoretical upper bound by balancing stability and acceleration. Particularly, the larger $p$ is, the faster the system evolves. Hence, the numerical integrator requires smaller step sizes to stabilize the process, but a smaller step size increases the number of iterations to achieve a target accuracy.  This tension is alleviated by Assumption \ref{assump:growth}. The larger $p$ is, the flatter the function $f$ is around its stationary points. In other words, Assumption~\ref{assump:growth} implies that as the iterates approach a minimum, the high order derivatives of the function $f$, in addition to the gradient, also converge to zero. Consequently, the trajectory slows down around the optimum and we can stably discretize the process with a large enough step size. This intuition ultimately translates into our main result.

\begin{theorem} \label{thm:main}
{\bf (Main Result)} Consider the second-order ODE in \eqref{eq:ode}. Suppose that the function $f$ is convex and Assumptions~\ref{assump:growth} and \ref{assump:differentiable} are satisfied. Further, let $s$ be the order of the Runge-Kutta integrator used in Algorithm~1, $N$ be the total number of iterations, and $x_0$ be the initial point. Also, let $\mathcal{E}_0:= f(x_0) - f(x^*) + \|x_0 - x^*\|^2+1$. Then, there exists a constant $C_1$ such that if we set the step size as $h=C_1N^{-1/(s+1)}(L+M+1)^{-1}\mathcal{E}_0^{-1}$, the iterate $x_N$ generated after running Algorithm~1 for $N$ iterations satisfies the inequality 
\begin{align}
f(x_N) - f(x^*) \le C_2 \mathcal{E}_0  \left[\tfrac{(L+M+1)\mathcal{E}_0}{N^{\frac{s}{s+1}}}\right]^p = \Oc\bigl(N^{-p\frac{s}{s+1}}\bigr),
\end{align}
where the constants $C_1$ and $C_2$ only depend on $s$, $p$, and the Runge-Kutta integrator. Since each iteration consumes $S$ gradient, $f(x_N) - f(x^*)$ will converge as $\Oc({S^{\frac{ps}{s+1}}N^{-\frac{ps}{s+1}}})$ with respect to the number of gradient evaluations. Note that for commonly used Runge-Kutta integrators, $S \le 8$.
\end{theorem}
The proof of this theorem is quite involved; we provide a sketch in Section~\ref{sec:proof}, deferring the detailed technical steps to the appendix. We do not need to know the constant $C_1$ exactly in order to set the step size $h$. Replacing $C_1$ by any smaller positive constant leads to the same polynomial rate.

Theorem~\ref{thm:main} indicates that if the objective has bounded high order derivatives and satisfies the flatness condition in Assumption~\ref{assump:growth} with $p>0$, then discretizing the ODE in \eqref{eq:ode} with a high order integrator results in an algorithm that converges to the optimal solution at a rate that is close to $\Oc({N^{-p}})$. In the following corollaries, we highlight two special instances of Theorem~\ref{thm:main}.

\begin{corollary}
If the function $f$ is convex with $L$-Lipschitz gradients and is $4^{\text{th}}$ order differentiable, then simulating the ODE~\eqref{eq:ode} for $p=2$ with a numerical integrator of order $s=2$ for N iterations results in the suboptimality bound
$$f(x_N) - f(x^*) \le \frac{C_2(f(x_0) - f(x^*) + \|x_0 - x^*\|^2+1)^3(L+M+1)^2}{N^{4/3}}.$$
\end{corollary}
Note that higher order differentiability allows one to use a higher order integrator, which leads to the optimal $\Oc({N^{-2}})$ rate in the limit. The next example is based on high order polynomial or $\ell_p$ norm.
\begin{corollary}
Consider the objective function $f(x) = \|Ax+b\|^4_4$. Simulating the ODE~\eqref{eq:ode} for $p=4$ with a numerical integrator of order $s=4$ for $N$ iterations results in the suboptimality bound
$$f(x_N) - f(x^*) \le \frac{C_2(f(x_0) - f(x^*) + \|x_0 - x^*\|^2+1)^5(L+M+1)^4}{N^{16/5}}.$$
\end{corollary}

\subsection{Logistic loss}\label{sec:logistic}
Discretizing logistic loss $f(x) = \log(1+e^{-w^T x})$ does not fit exactly into the setting of Theorem \ref{thm:main} due to nonexistence of $x^*$. This potentially causes two problems. First, Assumption~\ref{assump:growth} is not well defined. Second, the constant $\mathcal{E}_0$ in Theorem~\ref{thm:main} is not well defined. We explain in this section how we can modify our analysis to admit logistic loss by utilizing its structure of high order derivatives.

The first problem can be resolved by replacing $f(x^*)$ by $\inf_{x \in \reals^d}\! f(x)$ in Assumption~\ref{assump:growth}; then, the logistic loss satisfies Assumption \ref{assump:growth} with arbitrary integer $p > 0$. To approach the second problem, we replace $x^*$ by $\tilde{x}$ that satisfies the following relaxed inequalities. For some $\epsilon_1,\epsilon_2,\epsilon_3 < 1$ we have
 \begin{equation}
\inner{x-\tilde{x}, \nabla f(x)} \ge f(x) - f(\tilde{x}) - \epsilon_1,
\end{equation}
\begin{equation}
f(x) - f(\tilde{x}) \ge \tfrac{1}{L} \|\nabla^{(i)} f(x)\|^{\frac{p}{p-i}} - \epsilon_2, \qquad	f(\tilde{x}) -  \inf_{x \in \reals^d}\! f(x)  \leq \epsilon_3.
\end{equation}
As the inequalities are relaxed, there exists a vector $\tilde{x} \in \reals^d$ that satisfies the above conditions. If we follow the original proof and balance the additional error terms by picking $\tilde{x}$ carefully, we obtain
\begin{corollary}\label{cor:logistic}(Informal)
  If the  objective is $f(x) =\log(1+e^{-w^T x})$, then discretizing the ODE~ \eqref{eq:ode} with an order $s$ numerical integrator for $N$ iterations with step size $h = \Oc(N^{-1/(s+1)})$ results in a convergence rate of $\Oc({S^{p\frac{s}{s+1}}N^{-p\frac{s}{s+1}}})$.
\end{corollary}


\section{Proof of Theorem~\ref{thm:main}} \label{sec:proof}
We prove Theorem~\ref{thm:main} as follows. First(Proposition~\ref{prop_dec_energy}), we show that the suboptimality $f(x(t)) - f(x^*)$ along the continuous trajectory of the ODE~\eqref{eq:ode} converges to zero sufficiently fast. Second(Proposition~\ref{prop:discrete-error}), we bound the discretization error $\|\Phi_h(y_k) - \varphi_h(y_k)\|$, which measures the distance between the point generated by discretizing the ODE and the true continuous solution. Finally(Proposition~\ref{prop:main}), a bound on this error along with continuity of the Lyapunov function~\eqref{eq:lyapunov} implies that the suboptimality of the discretized sequence of points also converges to zero quickly.

Central to our proof is the choice of a Lyapunov function used to quantify progress. We propose in particular the Lyapunov function $\mathcal{E}: \reals^{2d+1}\to \reals_+$ defined as
\begin{equation}\label{eq:lyapunov}
  \mathcal{E}([v;x;t]) := \frac{t^2}{4p^2}\|v\|^2 + \Bigl\|x + \frac{t}{2p}v - x^*\Bigr\|^2 + t^p(f(x) - f(x^*)).
\end{equation}
The Lyapunov function~\eqref{eq:lyapunov} is similar to the ones used by~\citet{wibisono-variational,su-differential}, except for the extra term $\frac{t^2}{4p^2}\|v\|^2$. This term allows us to bound $\|v\|$ by $\Oc(\frac{\mathcal{E}}{t})$. This dependency is crucial for us to achieve the $O(N^{-2})$ bound(see Lemma~\ref{lemma:bound-high-order-F} for more details).

We begin our analysis with Proposition \ref{prop_dec_energy}, which shows that the function $\mathcal{E}$ is non-increasing with time, i.e., $\dot{\mathcal{E}}(y)\leq 0$. This monotonicity then implies that both $t^p(f(x) - f(x^*))$ and $ \frac{t^2}{4p^2}\|v\|^2$ are bounded above by some constants. The bound on $t^p(f(x) - f(x^*))$  provides a convergence rate of $\mathcal{O}(1/t^p)$ on the sub-optimality $f(x(t)) - f(x^*)$. It further leads to an upper-bound on the derivatives of the function $f(x)$ in conjunction with Assumption~\ref{assump:growth}. 

\begin{proposition}[Monotonicity of $\cal E$]\label{prop_dec_energy}
Consider the vector $y = [v;x;t]\in \mathbb{R}^{2d+1}$ as a trajectory of the dynamical system~\eqref{eq:dynamics}. Let the Lyapunov function $\mathcal{E}$ be defined by~\eqref{eq:lyapunov}. Then, for any trajectory $y = [v;x;t]$, the time derivative $\dot{\mathcal{E}}(y) $ is non-positive and bounded above; more precisely, 
\begin{equation}
\dot{\mathcal{E}}(y) \le - \frac{t}{p}\|v\|^2.
\end{equation}
\end{proposition}
The proof of this proposition follows from convexity and \eqref{eq:ode}; we defer the details to Appendix~\ref{supp_prop_dec_energy}.


Next, to bound the Lyapunov function for numerical solutions, we need to bound the distance between points in the discretized and continuous trajectories. As in Section~\ref{sec:RK}, for the dynamical system $\dot{y} = F(y)$, let $\Phi_h(y_0)$ denote the solution generated by a numerical integrator starting at point $y_0$ with step size $h$. Similarly, let $\varphi_h(y_0)$ be the corresponding true solution to the ODE. An ideal numerical integrator would satisfy $\Phi_h(y_0)= \varphi_h(y_0)$; however, due to discretization error there is always a difference between $\Phi_h(y_0)$ and $\varphi_h(y_0)$ determined by the order of the integrator as in~\eqref{eq:disc-error}. Let $\{y_k\}_{i=0}^N$ be the sequence of points generated by the numerical integrator, that is, $y_{k+1} = \Phi_h(y_k)$. In the following proposition, we derive an upper bound on the resulting discretization error $\|\Phi_h(y_k) - \varphi_h(y_k)\|$.

\begin{proposition}[Discretization error]
  \label{prop:discrete-error}
Let $y_k = [v_k;x_k; t_k]$ be the current state of the dynamical system $\dot{y} = F(y)$ defined in~\eqref{eq:dynamics}. Suppose $x_k \in \mathcal{S}$ defined in ~\eqref{eq:set-S}. If we use a Runge-Kutta integrator of order $s$ to discretize the ODE for a single step with a step size $h$ such that $h \le  \min\{0.2,\frac{1}{(1+\kappa)C(1+\mathcal{E}(y_k))(M+L+1)}\}$, then 
\begin{equation}\label{claim_dis_error}
\|\Phi_h(y_k) - \varphi_h(y_k)\| \le C' h^{s+1}(M\!+\!L\!+\!1)\left[ \frac{[(1+\mathcal{E}(y_k))]^{s+1}}{t_k} + h \frac{[(1+\mathcal{E}(y_k))]^{s+2}}{t_k}\right],
\end{equation}
where the constants $C$, $\kappa$, and $C'$ only depend on $p$, $s$, and the integrator.
\end{proposition}
The proof of Proposition~\ref{prop:discrete-error} is the most challenging part in proving Theorem~\ref{thm:main}. Details may be found in Appendix~\ref{ap:discretization}.
The key step is to bound $\|\frac{\partial^{s+1}}{\partial h^{s+1}}[\Phi_h(y_k) - \varphi_h(y_k)]  \|$. To do so, we first bound the high order derivative tensor $\|\nabla^{(i)}f\|$ using Assumption~\ref{assump:growth} and Proposition~\ref{prop_dec_energy} within a region of radius $R$. By carefully selecting  $R$, we can show that for a reasonably small $h$, $\Phi_h(y_k)$ and $\varphi_h(y_k)$ is constrained in the region. Second, we need to compute the high order derivatives of $\dot{y} = F(y)$ as a function of $\nabla^{(i)}f$ which is bounded in the region of radius R. As shown in Appendix \ref{sec:elementary}, the expressions for higher derivatives become quite complicated as the order increases.  We approach this complexity by using the notation for elementary differentials (see Appendix \ref{sec:elementary}) adopted from \citep{hairer-textbook}; we then induct on the order of the derivatives to bound the higher order derivatives. The flatness assumption (Assumption~\ref{assump:growth}) provides bounds on the operator norm of high order derivatives relative to the objective function suboptimality, and hence proves crucial in completing the inductive step.

By the conclusion in Proposition~\ref{prop:discrete-error} and continuity of the Lyapunov function $\mathcal{E}$, we conclude that the value of $\mathcal{E}$ at a discretized point is close to its continuous counterpart. Using this observation, we expect that the Lyapunov function values for the points generated by the discretized ODE do not increase significantly. We formally prove this key claim in the following proposition.

\begin{proposition}\label{prop:main}
Consider the dynamical system $\dot{y} = F(y)$ defined in~\eqref{eq:dynamics} and the Lyapunov function~$\mathcal{E}$ defined in \eqref{eq:lyapunov}. Let $y_0$ be the initial state of the dynamical system and $y_N$ be the final point generated by a Runge-Kutta integrator of order $s$ after $N$ iterations. Further, suppose that  Assumptions~\ref{assump:growth} and \ref{assump:differentiable} are satisfied. Then, there exists a constant $\tilde{C}$ determined by $p, s$ and the numerical integrator, such that if the step size $h$ satsfies $h = \tilde{C} \frac{N^{-1/(s+1)}}{(L+M+1)(e\mathcal{E}(y_0)+1)}$, then we have
\begin{equation}\label{energy_bound_claim}
\mathcal{E}(y_N) \le \exp(1)\ \mathcal{E}(y_0) + 1.
\end{equation}
\end{proposition}
Please see Appendix~\ref{ap:bound-discrete-lyapunov} for a proof of this claim.

Proposition~\ref{prop:main} shows that the value of the Lyapunov function $\mathcal{E}$ at the point $y_N$ is bounded above by a constant that depends on the initial value $\mathcal{E}(y_0)$. Hence, if the step size $h$ satisfies the required condition in Proposition~\ref{prop:main}, we can see that
\begin{align}\label{disc_sub_opt_result_1}
  f(x_N) - f(x^*) \le \tfrac{\mathcal{E}(y_N)}{t_N^p} \le \tfrac{e\mathcal{E}(y_0) + 1}{(1+Nh)^p}.
\end{align}
The first inequality in~\eqref{disc_sub_opt_result_1} follows from the definition of the $\cal E$~\eqref{eq:lyapunov}. Replacing the step size $h$ in \eqref{disc_sub_opt_result_1} by the choice used in Proposition~\ref{prop:main} yields
\begin{align}\label{disc_sub_opt_result_2}
f(x_N) - f(x^*) \le \frac{(L+M+1)^p(e\mathcal{E}(y_0) + 1)^{p+1}}{\tilde{C}N^{p\frac{s}{s+1}}},
\end{align}
and the claim of Theorem~\ref{thm:main} follows.

\emph{Note:} The dependency of the step size $h$ on the degree of the integrator $s$ suggests that an integrator of higher order allows for larger step size and therefore faster convergence rate.


\section{Numerical experiments} \label{sec:exp}
In this section, we perform a series of numerical experiments to  study the performance of the proposed scheme for minimizing convex functions through the direct discretization (DD) of the ODE in \eqref{eq:ode} and compare it with gradient descent (GD) as well as Nesterov's accelerated gradient (NAG). All figures in this section are on log-log scale. For each method tested, we empirically choose the largest step size among $\{10^{-k}| k \in \mathcal{Z} \}$ subject to the algorithm remaining stable in the first 1000 iterations.

\subsection{Quadratic functions} \label{sec:exp-quad}

We now verify our theoretical results by minimizing a quadratic convex function of the form $f(x) = \|Ax-b\|^2$ by simulating the ODE in \eqref{eq:ode} for the case that $p=2$, i.e., 
$$\ddot{x}(t) + \frac{5}{t}\dot{x}(t) + 4 \nabla f(x(t)) = 0,$$
where $A\in \mathbb{R}^{10\times 10}$ and  $b\in  \mathbb{R}^{10}$. The first five entries of $b=[b_1;\dots;b_{10}]$ are valued $0$ and the rest are $1$. Rows $A_i$ in $A$ are generated by an $i.i.d$ multivariate Gaussian distribution conditioned on $b_i$. The data is linearly separable.  Note that the quadratic objective $f(x) = \|Ax-b\|^2$ satisfies the condition in Assumption \ref{assump:growth} with $p=2$. It is also clear that it satisfies the condition in Assumption~\ref{assump:differentiable} regarding the bounds on higher order derivatives.

\begin{figure}[t]
    \centering
    \begin{subfigure}[b]{0.5\textwidth}
	\centering
	\includegraphics[width=\textwidth]{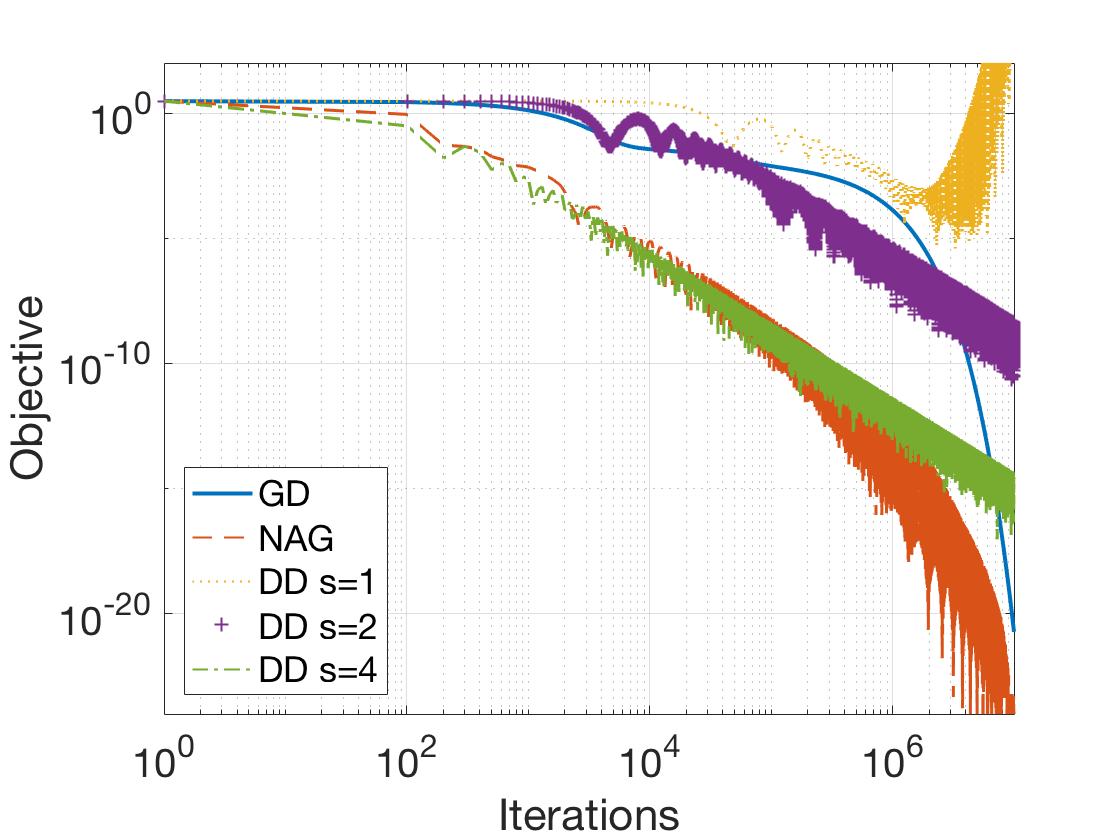}
	\caption{Quadratic objective}
	\end{subfigure}%
	~ 
	\begin{subfigure}[b]{0.5\textwidth}
	\centering
	\includegraphics[width=\textwidth]{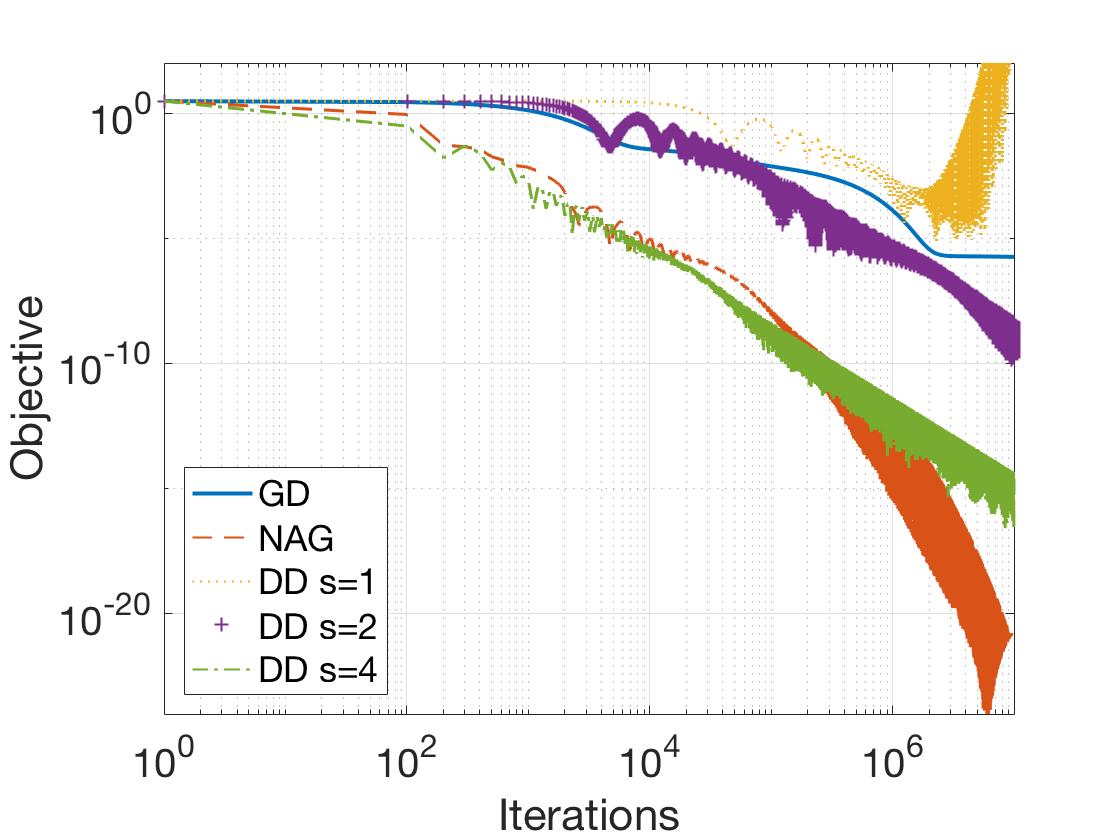}
	\caption{Objective as in \eqref{eq:24-objective}}
	\end{subfigure}

    \caption{Convergence paths of \gd, \nag, and the proposed simulated dynamical system with integrators of degree $s=1$, $s=2$, and $s=4$. The objectives satisfy Assumption~\ref{assump:growth} with p=2. }
    \label{fig:quadratic}
\end{figure}

Convergence paths of GD, NAG, and the proposed DD procedure for minimizing the quadratic function $f(x) = \|Ax-b\|^2$ are demonstrated in Figure~\ref{fig:quadratic}(a). For the proposed method we consider integrators with different degrees, i.e., $s\in\{1,2,4\}$.  Observe that GD eventually attains linear rate since the function is strongly convex around the optimal solution. \nag displays local acceleration close to the optimal point as mentioned in \citep{su-differential,attouch-local-nag}. For DD, if we simulate the ODE with an integrator of order $s=1$, the algorithm is eventually unstable. This result is consistent with the claim in \citep{wibisono-variational} and our theorem that requires the step size to scale with $\Oc({N^{-0.5}})$. Notice that using a higher order integrator leads to a stable algorithm. Our theoretical results suggest that the convergence rate for $s\in\{1,2,4\}$ should be worse than $\Oc({N^{-2}})$ and one can approach such rate by making $s$ sufficiently large. However, as shown in Figure~\ref{fig:quadratic}(a), in practice with an integrator of order $s=4$, the DD algorithm achieves a convergence rate of $\Oc({N^{-2}})$. Hence, our theoretical convergence rate in Theorem 1 might be conservative.

We also compare the performances of these algorithms when they are used to minimize 
\begin{equation}\label{eq:24-objective}
f([x_1, x_2]) = \|Ax_1-b\|^2 + \|C x_2 - d\|_4^4.
\end{equation}
Matrix $C$ and vector $d$ are generated similarly as $A$ and $b$. The result is shown in Figure \ref{fig:quadratic}(b). As expected, we note that \gd no longer converges linearly, but the other methods converge at the same rate as in Figure \ref{fig:quadratic}(a).

\subsection{Decoupling ODE coefficients with the objective} \label{sec:decouple}

\begin{figure}[t]
    \centering
    \includegraphics[width=0.5\textwidth]{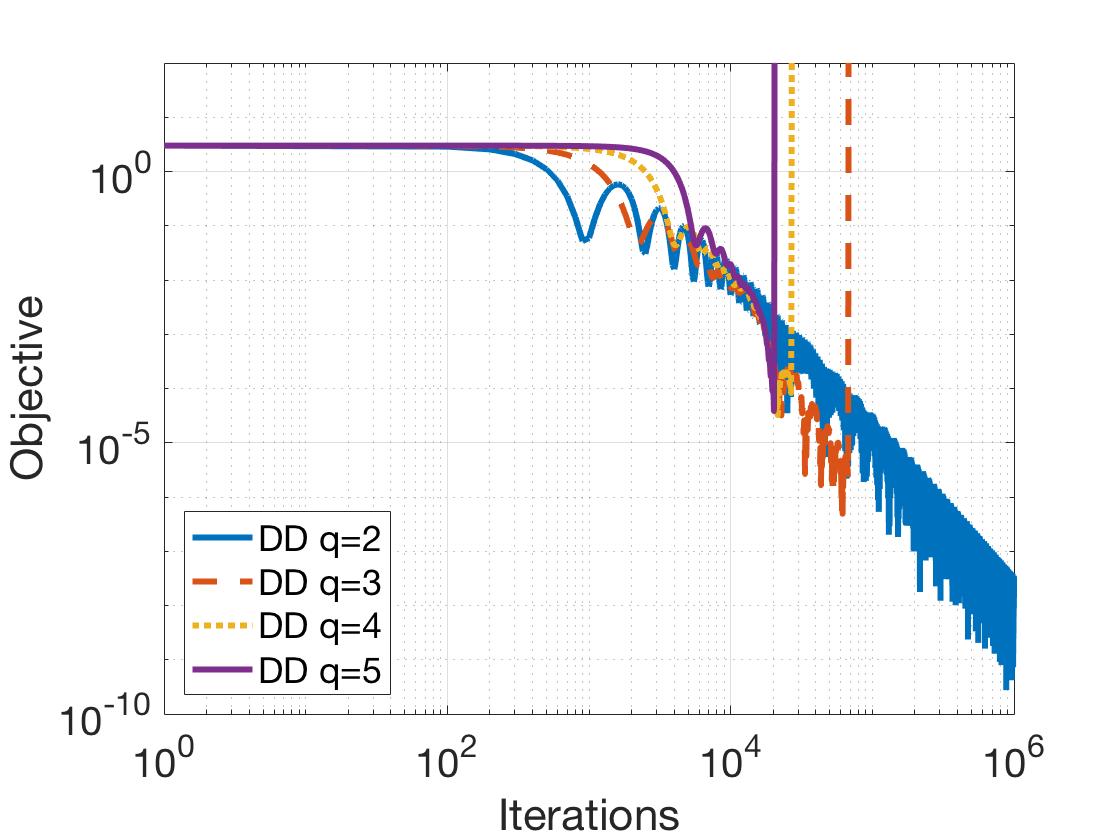}
    \caption{Minimizing quadratic objective by simulating different ODEs with the RK44 integrator ($4^{th}$ order). In the case when $p=2$, the optimal choice for q is 2.}
    \label{fig:multiple-q}
\end{figure}

Throughout this paper, we assumed that the constant $p$ in (\ref{eq:ode}) is the same as the one in Assumption~\ref{assump:growth} to attain the best theoretical upper bounds. In this experiment, however, we empirically explore the convergence rate of discretizing the ODE 
$$\ddot{x}(t) + \frac{2q+1}{t}\dot{x}(t) + q^2 t^{q-2} \nabla f(x(t)) = 0,$$
when $q \neq p$. In particular, we use the same quadratic objective $f(x) = \|Ax-b\|^2$
as in the previous section. This objective satisfies Assumption \ref{assump:growth} with $p=2$. We simulate the ODE with different values of $q$ using the same numerical integrator with the same step size. Figure~\ref{fig:multiple-q} summarizes the experimental results. We observe that when $q > 2$, the algorithm diverges. Even though the suboptimality along the continuous trajectory will converge at a rate of $\Oc({t^{-p}})=\Oc({t^{-2}})$, the discretized sequence cannot achieve the lower bound which is of $\Oc({N^{-2}})$.

\subsection{Beyond Nesterov's acceleration}
\begin{figure}[t]
    \centering

	\begin{subfigure}[b]{0.5\textwidth}
		\centering
		\includegraphics[width=\textwidth]{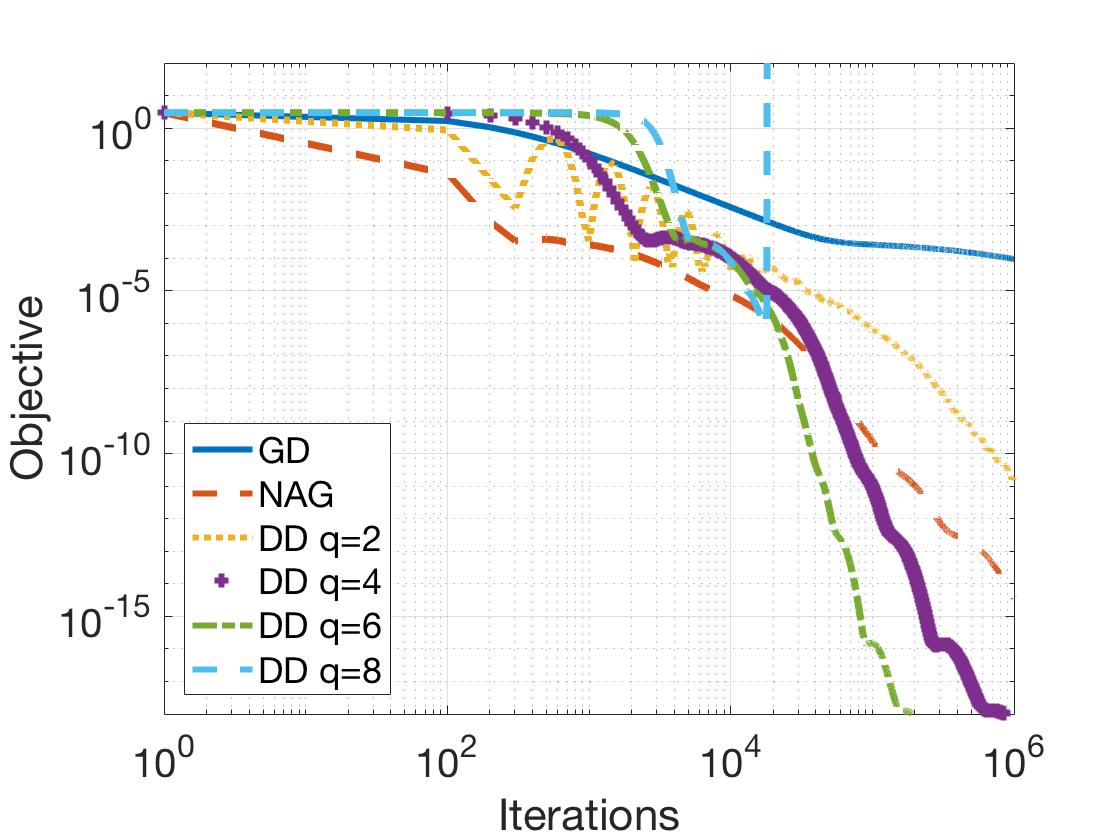}
		\caption{The objective is an $\ell_4$ norm.}
	\end{subfigure}%
	~ 
	\begin{subfigure}[b]{0.5\textwidth}
		\centering
		\includegraphics[width=\textwidth]{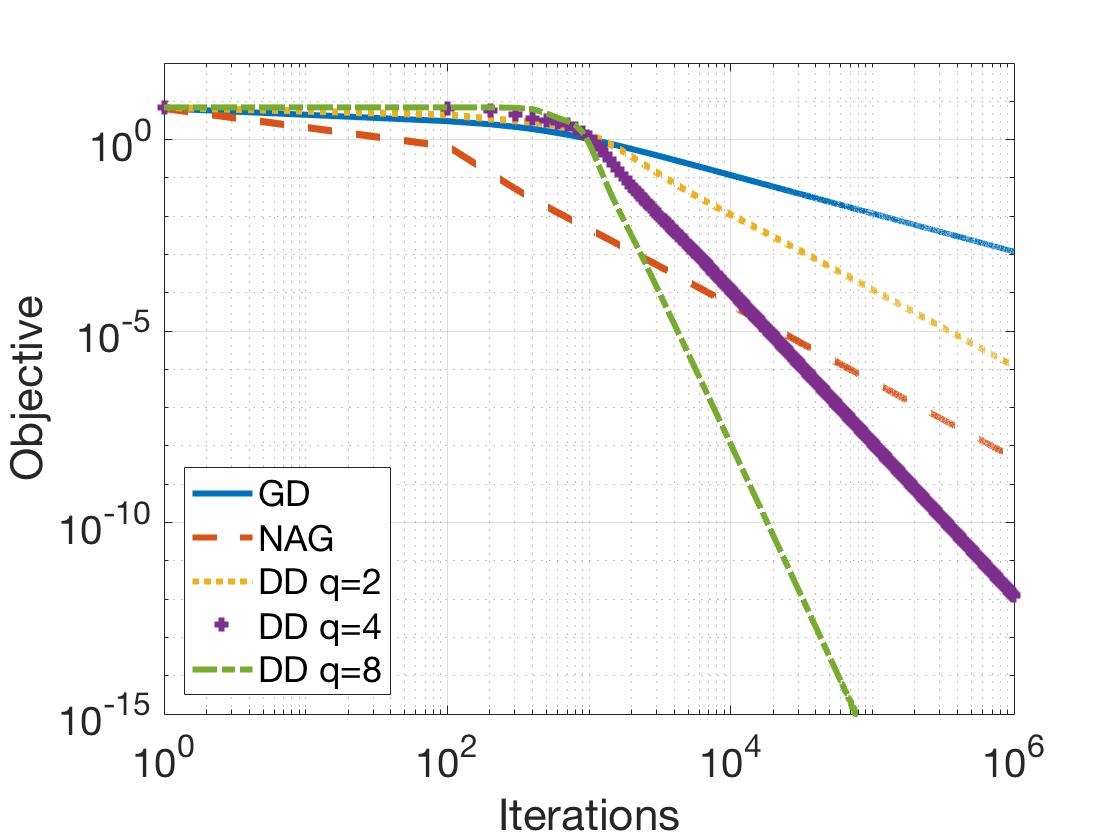}
		\caption{The objective is a logistic loss.}
	\end{subfigure}
	
    \caption{Experiment results for the cases that Assumption \ref{assump:growth} holds for $p>2$.  }
    \label{fig:beyond}
\end{figure}
In this section, we discretize ODEs with objective functions that satisfy Assumption \ref{assump:growth} with $p>2$. For all ODE discretization algorithms, we use an order-$2$ RK integrator that calls the gradient oracle twice per iteration. Then we run all algorithms for $10^6$ iterations and show the results in Figure \ref{fig:beyond}. 
As shown in Example \ref{ex:4th}, the objective
\begin{equation}
f(x) = \|Ax-b\|^4_4
\end{equation}
satisfies Assumption $\ref{assump:growth}$ for $p=4$. By Theorem \ref{thm:main} if we set $q=4$, we can achieve a convergence rate close to the rate $\Oc({N^{-4}})$. We run the experiments with different values of $q$ and summarize the results in Figure~\ref{fig:beyond}(a). Note that when $q>2$, the convergence of direct discretization methods is faster than \nag. Interestingly, when $q=6>p=4$, the discretization is still stable with convergence rate roughly $\Oc({N^{-5}})$. This suggests that our theorem may be conservative.

We then simulate the ODE for the objective function 
$$f(x) = \sum_{i=1}^{10} \log(1+e^{-w_i^Tx}),$$
for a dataset of linearly separable points. The data points are generated in the same way as in Section~\ref{sec:exp-quad}. As shown in Section~\ref{sec:logistic}, it satisfies Assumption~\ref{assump:growth} for arbitrary $p>0$. As shown in Figure \ref{fig:beyond}(b), the objective decreases faster for larger $q$; this verifies Corollary~\ref{cor:logistic}.


\section{Discussion}
This paper specifies sufficient conditions for stably discretizing an ODE to obtain accelerated first-order (i.e., purely gradient based) methods. Our analysis allows for the design of optimization methods via direct discretization using Runge-Kutta integrators based on the flatness  of  objective functions. complementing the existing studies that derive ODEs from optimization methods, we show that one can prove convergence rates of a optimization algorithms by leveraging properties of its ODE representation. We hope that this perspective will lead to more general results.

In addition, we identified a new condition in Assumption \ref{assump:growth} that quantifies the \emph{local flatness} of convex functions. At first, this condition may appear counterintuitive, because gradient descent actually converges fast when the objective is \emph{not} flat and the progress slows down if the gradient vanishes close to the minimum. However, when we discretize the ODE, the trajectories with vanishing gradients oscillate slowly, and hence allow stable discretization with large step sizes, which ultimately allows us to achieve acceleration. We think this high-level idea, possibly as embodied by Assumption~\ref{assump:growth} could be more broadly used in analyzing and designing other optimization methods.

Based on the above two points, this paper contains both positive and negative message for the recent trend in ODE interpretation of optimization methods. On one hand, it shows that with careful analysis, discretizing ODE can preserve some of its trajectories properties. On the other hand, our proof suggests that nontrivial additional conditions might be required to ensure stable discretization. Hence, designing an ODE with nice properties in the continuous domain doesn't guarantee the existence of a practical optimization algorithm.

Although our paper answers a fundamental question regarding the possibility of obtaining accelerated gradient methods by directly discretizing second order ODEs (instead of reverse engineering Nesterov-like constructions), it does not fully explain  acceleration. First, unlike Nesterov's accelerated gradient method that only requires first order differentiability, our results require the objective function to be $(s+2)$-times differentiable (where $s$ is the order of the integrator). Indeed, the precision of numerical integrators only increases with their order when the function is sufficiently differentiable. This property inherently limits our analysis. Second, while we achieve the $\Oc(N^{-2})$ convergence rate, some of the constants in our bound are loose (e.g., for squared loss and logistic regression they are quadratic in $L$ versus linear in $L$ for \nag). Achieving the optimal dependence on initial errors $f(x_0)-f(x^*)$, the diameter $\|x_0-x^*\|$, as well as constants $L$ and $M$ requires further investigation.


\section*{Acknowledgement}
AJ and SS acknowledge support in part from DARPA FunLoL, DARPA Lagrange; AJ also acknowledges support from an ONR Basic Research Challenge Program, and SS acknowledges support from NSF-IIS-1409802.

{
  \bibliography{bibliography}
  \setlength{\bibsep}{1pt}
\bibliographystyle{abbrvnat}
}
\clearpage

\appendix

\section{Proof of Proposition~\ref{prop_dec_energy}}\label{supp_prop_dec_energy}

According to the dynamical system in \eqref{eq:dynamics} we can write
\begin{align}
\dot{x} = v ,\qquad 
\ddot{x} = \dot{v} = -\frac{2p+1}{t}v - p^2t^{p-2}\nabla f(x).
\end{align}
Using these definitions we can show that 
\begin{align}
\dot{\mathcal{E}} = &\frac{t^2}{4p^2}\inner{2v, \dot{v}} + \frac{2t}{4p^2}\inner{v, v} + 2\inner{x + \frac{t}{2p}v - x^*, \dot{x} + \frac{\dot{x}}{2p} + \frac{t}{2p}\ddot{x}}+ t^p\inner{\nabla f(x), \dot{x}}
\nonumber\\
&\quad  + pt^{p-1}(f(x) -  f(x^*)) \nonumber\\
=& \frac{2t^2}{4p^2} \inner{\dot{x}, \ddot{x} + \frac{2p+1}{t}\dot{x}} - \frac{2t^2}{4p^2} \inner{\dot{x}, \frac{2p}{t}\dot{x}}+ 2\frac{t}{2p}\inner{x + \frac{t}{2p}\dot{x} - x^*, \ddot{x} + \frac{2p+1}{t}\dot{x}} \nonumber\\
&\quad + t^p\inner{\nabla f(x), \dot{x}} + pt^{p-1}(f(x) -  f(x^*)) \nonumber\\
=& \frac{t^2}{2p^2}\inner{\dot{x}, -p^2t^{p-2} \nabla f} - \frac{t}{p}\|\dot{x}\|^2 + \frac{t}{p}\inner{x + \frac{t}{2p}\dot{x} - x^*, -p^2t^{p-2} \nabla f} \nonumber \\
& \quad+ t^p\inner{\nabla f(x), \dot{x}} + pt^{p-1}(f(x) -  f(x^*))\nonumber \\
=& - \frac{t}{p}\|\dot{x}\|^2  + pt^{p-1}(f(x) -  f(x^*)) - pt^{p-1}\inner{x - x^*, \nabla f} \nonumber\\
\le & - \frac{t}{p}\|\dot{x}\|^2.
\end{align}
The equalities follows from rearrangement and \eqref{eq:ode}. The last inequality holds due to convexity.

\section{Proof of Proposition~\ref{prop:discrete-error} (Discretization Error)}\label{ap:discretization}
In this section, we aim to bound the difference between the true solution defined by the ODE and the point generated by the integrator, i.e., $\|\Phi_h(y_c) - \varphi_h(y_c)\|$. Since the integrator has order $s$, the difference $\Delta(h):= \|\Phi_h(y_c) - \varphi_h(y_c)\|$ should be proportional to $h^{s+1}$. Here, we intend to formally derive an upper bound of $\mathcal{O}(h^{s+1})$ on $\Delta(h)$.


We start by introducing some notations. Given a vector $y = [v;x;t]\in \reals^{2d+1}$, we define the following projection operators
\begin{align}
\pi_x(y) = x \in \reals^d, \quad \pi_v(y) = v \in \reals^d, \quad \pi_t(y) = t \in \reals, \quad \pi_{v,x}(y) = \vektor{v,x,} \in \reals^{2d}.
\end{align}
We also define the set $B(x_c, R) $ which is a ball with center $x_c$ and radius $R$ as 
\begin{equation}\label{eq:B_R}
	B(x_c, R) = \{x \in \reals^d | \|x-x_c\| \le R \} ,
\end{equation}
and define the set $U_{R, 0.2}(y_c)$ as
\begin{equation} \label{eq:U}
U_{R, 0.2}(y_c) = \{y=[v; x; t] | \|v-v_c\|\le R, \|x - x_c\| \le R, |t-t_c|\le 0.2\}.
\end{equation}

In the following Lemma, we show that if we start from the point $y_c$ and choose a sufficiently small stepsize, the true solution defined by the ODE $\varphi_h(y_0)$ and the point generated by the integrator $\Phi_h(y_c)$ remain in the set $U_{R, 0.2}(y_c)$.

\begin{lemma} \label{lemma:stay_in_the_set}
Let $ y \in U_{R, 0.2}(y_c)$ where $y_c = [v_c; x_c; t_c]$, $t_c \ge 1$, and $R = \frac{1}{t_c}$. Suppose that $ B(x_c, R) \subseteq \mathcal{A}$ (defined in~\eqref{eq:set-A}) and hence Assumptions \ref{assump:growth} and \ref{assump:differentiable} are satisfied.  If $h \le \min\{0.2,\frac{1}{(1+\kappa)C(\mathcal{E}(y_c) + 1)(L+M+1)}\}$, the true solution defined by the ODE $\varphi_h(y_0)$ and the point generated by the integrator $\Phi_h(y_c)$ remain in the set $U_{R, 0.2}(y_c)$, i.e., 
\begin{align}
\varphi_h(y_c) \in U_{R, 0.2}(y_c), \qquad \Phi_h(y_c) \in U_{R, 0.2}(y_c),
\end{align}
where $\kappa$ is a constant determined by the Runge-Kutta integrator. In addition, the intermediate points $g_i$ defined in Definition~\ref{def:rk} also belong to the set $U_{R, 0.2}(y_c)$.
\end{lemma}

\begin{proof}

Note that $\forall y \in \reals^{2d+1}$, $	\|\pi_t F(y)\| = 1$. Clearly when $h \le 0.2$,
\begin{equation}
\pi_t\varphi_h(y_c)-y_c = h \le 0.2.
\end{equation}
Similarly, for any integrator that is at least order $1$,
\begin{equation}
\pi_t\Phi_h(y_c)-y_c = h \le 0.2.
\end{equation}
Therefore, we only need to focus on bounding the remaining coordinates.

By Lemma~\ref{lemma:bound-F}, we have that when $y \in  U_{R, 0.2}(y_c)$,
\begin{equation}
	\|\pi_{v,x} F(y)\| \le \frac{C(\mathcal{E}(y_c) + 1)(L+M+1)}{t_c}.
\end{equation}

By definition \ref{def:rk},
\begin{equation*}
g_i = y_k + h \sum_{j=1}^{i-1} a_{ij} F(g_j) \qquad \Phi_h(y_k) = y_k + h \sum_{i=0}^{s-1} b_i F(g_i).
\end{equation*}
Let $\kappa = \max\{\sum_j |a_{ij}|, \sum |b_i|\}$, we have that when $h \le \min\{0.2,R/[\kappa\frac{C(\mathcal{E}(y_c) + 1)(L+M)}{t_c}]\}$,
\begin{equation}
g_i \in U_{R, 0.2}(y_c)\qquad \Phi_h(y_c) \in U_{R, 0.2}(y_c).
\end{equation}

By fundamental theorem of calculus, we have that
\begin{equation}
	\varphi_h(y_c) = y_c + \int_0^h F(\varphi_t(y_c)) dt \in U_{R, 0.2}(y_c).
\end{equation}
Rearrange and apply Cauchy-Schwarz, we get
\begin{equation}
\|\pi_{v,x}[\varphi_h(y_c) - y_c]\| \le \int_0^h \|\pi_{v,x} F(\varphi_t(y_c))\| dt \in U_{R, 0.2}(y_c).
\end{equation}
By mean value theorem and proof of contradiction, we can show that when $h \le \min\{0.2,R/\frac{C(\mathcal{E}(y_c) + 1)(L+M)}{t_c}\}$,
\begin{equation}
 \int_0^h \|\pi_{v,x}F(\varphi_t(y_c))\| dt \le R.
\end{equation}
In particular, if $ \int_0^h \|\pi_{v,x}F(\varphi_t(y_c))\| dt \ge R$, then exists $y_1$ and $h_0 < h$ such that $\|y_1-y_c\|=R$ and  $y_1= y_c + \int_0^{h_0} F(\varphi_t(y_c)) dt $. By mean value theorem, this implies that exist $y \in U_{R, 0.2}(y_c)$ such that $\|\pi_{v,x} F(y)\| > \frac{C(\mathcal{E}(y_c) + 1)(L+M+1)}{t_c}$, which contradicts Lemma~\ref{lemma:bound-F}.

Therefore we proved that
\begin{equation}
\varphi_h(y_c) \in U_{R, 0.2}(y_c).
\end{equation}
\end{proof}


The result in Lemma~\ref{lemma:stay_in_the_set} shows that $\varphi_h(y_c)$ and $\Phi_h(y_c)$ remain in the set $U_{R, 0.2}(y_c)$. In addition, we can bound the operator norm of $\nabla^{(i)}f$ in $B(x_c, R)$ by Lemma~\ref{lemma:bound-high-order-f}.  Since $\frac{\partial^q \varphi_h(y_c)}{\partial h^q}$ is a function of $\nabla^{(i)}f$, we can show in Lemma \ref{lemma:bound-high-order-F} that the $(s+1)_{th}$ derivative of $\varphi_h(y_c)$ and $\Phi_h(y_c)$ are bounded above by 
\begin{equation}\label{eq:bound_on_s_der_a}
\left\| \frac{\partial^q \varphi_h(y_c)}{\partial h^q} \right\| \le  \frac{C_0 [\mathcal{E}(y_c) + 1]^{q}(L+M+1)^q}{t_c},
\end{equation}
and
\begin{equation}\label{eq:bound_on_s_der_b}
\left\| \frac{\partial^q \Phi_h(y_c)}{\partial h^q} \right\| \le  \frac{C_1 [1 + \mathcal{E}(y_c)]^q (L+M+1)^q + C_2 h[1+\mathcal{E}(y_c)]^{q+1}(L+M+1)^{p+1}}{t_c}. 
\end{equation}

Since the integrator has order $s$, we can write 
\begin{equation}
\frac{\partial^i}{\partial h^i}[\Phi_h(y_k) - \varphi_h(y_k)] = 0\ \ \ \text{for} \ i = 1, ..., s.
\end{equation}
Therefore, the difference between the true solution $\varphi_h(y_c)$ defined by the ODE  and the point $\Phi_h(y_c)$ generated by the integrator can be upper bounded by
\begin{equation}\label{close}
 \|\Phi_h(y_c) - \varphi_h(y_c)\| \leq \left(\left\| \frac{\partial^{s+1} \varphi_h(y_k)}{\partial h^{s+1}} \right\| + \left\|\frac{\partial^{s+1} \Phi_h(y_k)}{\partial h^{s+1}} \right\|\right) h^{s+1}
\end{equation}
Replacing the norms on the right hand side of \eqref{close} by their upper bounds in \eqref{eq:bound_on_s_der_a} and \eqref{eq:bound_on_s_der_b} implies that
\begin{align}\label{final_step_bound}
 \|\Phi_h(y_c) - \varphi_h(y_c)\| & \le h^{s+1}\left[  \frac{(C_{0}+C_1) [\mathcal{E}(y_c) + 1]^{s+1}(M+L+1)^{s+1}}{t_c} \right]  
\nonumber\\
&\qquad  + h^{s+2}\left[ \frac{ C_{2}[1+\mathcal{E}(y_c)]^{s+2}(M+L+1)^{s+2}}{t_c} \right]. 
\end{align}
By replacing $y_c=[v_c;x_c;t_c]$ in \eqref{final_step_bound} by $y_k=[v_k;x_k;t_k]$ the claim in \eqref{claim_dis_error} follows.

\section{Proof of Proposition \ref{prop:main} (Analysis of discrete Lyapunov functions)} \label{ap:bound-discrete-lyapunov}
As defined earlier in Section~\ref{sec:proof}, $\Phi_h(y_k)$ is the solution generated by the numerical integrator, and $\varphi_h(y_k)$ is a point on the trajectory of the ODE. $y_c = [\vec{0}; x_c; 1]$ is the initial point of the ODE. Recall that $\{y_k\}_{i=0}^N$ is the sequence of points produced by the numerical integrator, i.e., $y_{k+1} = \Phi_h(y_k)$. 

To simplify the notation, we let $E_k = \mathcal{E}(y_k)$, $E_{k+1} = \mathcal{E}(\Phi_h(y_k))$, $\tilde{y} = \varphi_h(y_k) = [\tilde{v}; \tilde{x}; t+h] $, $\hat{y} = \Phi_h(y_k) = [\hat{v}; \hat{x}; t+h] $. 

We want to prove by induction on $k=0,1,...,N$ that
\begin{equation}\label{claim_of_induction}
E_k \le (1+\frac{1}{N})^k E_0 + \frac{k}{N}.
\end{equation}
The base case $E_0 \le E_0$ is trivial. Now let's assume by induction that the inequality in \eqref{claim_of_induction} holds for $k=j$, i.e., 
\begin{equation}\label{ass_induc}
E_j \le (1+\frac{1}{N})^j E_0 + \frac{j}{N}.
\end{equation}
By this assumption, we know that  $f(x_k) \le \frac{eE_0 + 1}{t_k^p} \le eE_0 + 1$ and hence $x_k \in \mathcal{S}$ defined in ~\eqref{eq:set-S}.
Note that $R=\frac{1}{t_k}\le 1$. We then have
\begin{equation}\label{eq:in-set-A}
B(x_k, R) \subseteq B(x_k, 1) \in \mathcal{A}
\end{equation}
 for $\mathcal{A}$ defined in \eqref{eq:set-A}.
By assumption in Proposition~\ref{prop_dec_energy},
\begin{align} \label{eq:h-ineq}
h \le  0.2,\qquad h \le \frac{1}{(1+\kappa)C(eE_0 + 2)(L+M+1)}.
\end{align}

By utilizing the bound on $ \|\Phi_h(y_k) - \varphi_h(y_k)\|$ and the continuity of $\mathcal{E}(y)$, we show in Lemma~\ref{lemma:bound-discrete-lyapunov} that the discretization error of  $\|\mathcal{E}(\hat{y}) - \mathcal{E}(\tilde{y})\|$ is upper bounded by 
\begin{align}\label{inter_result}
&\|\mathcal{E}(\Phi_h(y_k)) - \mathcal{E}(\varphi_h(y_k))\| \\
&\le C' h^{s+1}[(1+E_k)^{s+1}(L+M+1)^{s+1} \! + h(1+E_k)^{s+2}(L+M+1)^{s+2}](E_k + E_{k+1} + 1),\nonumber
\end{align}
under conditions in \eqref{eq:in-set-A} and \eqref{eq:h-ineq}. $C'$ only depends on $p, s$ and the numerical integrator.

We proceed to prove the inductive step. Start by writing $E_{k+1}=\mathcal{E}(\Phi_h(y_k))$ as 
\begin{equation}
\mathcal{E}(\Phi_h(y_k)) = \mathcal{E}(y_{k}) +  \mathcal{E}(\varphi_h(y_k)) - \mathcal{E}(y_{k}) +  \mathcal{E}(\Phi_h(y_k)) -  \mathcal{E}(\varphi_h(y_k)).
\end{equation}
According to Proposition~\ref{prop_dec_energy}, $ \mathcal{E}(\varphi_h(y_k)) - \mathcal{E}(y_{k}) \le 0$. Therefore,
\begin{equation}
E_{k+1} \le E_k + \|\mathcal{E}(\hat{y}) - \mathcal{E}(\tilde{y})\|.
\end{equation}
Replace the norm $\|\mathcal{E}(\hat{y}) - \mathcal{E}(\tilde{y})\|=\|\mathcal{E}(\Phi_h(y_k)) - \mathcal{E}(\varphi_h(y_k))\|$ by its upper bound \eqref{inter_result} to obtain
\begin{equation}\label{eq_for_induction}
E_{j+1} \le E_j + C h^{s+1}[(1+E_j)^{s+1}(L+M+1)^{s+1} \! + h(1+E_j)^{s+2}(L+M+1)^{s+2}](E_j + E_{j+1} + 1).
\end{equation}
Before proving the inductive step, we need to ensure that the step size $h$ is sufficiently small. Here, we further add two more $j$-independent conditions on the choice of step size $h$. In particular, we assume that 
\begin{align} \label{eq:h-ineq_2}
h \le \frac{1}{eE_0 + 2},\qquad  h^{s+1} \le \frac{1}{3(1+C^{-1})C'N(eE_0 + 2)^{s+1}(L+M+1)^{s+1} }.
\end{align}

Note that since we want show the claim in \eqref{claim_of_induction} for $k=1,\dots,N$, in inductive assumptions we have that $j\leq N-1$. Now we proceed to show that if the inequality in \eqref{claim_of_induction} holds for $k=j$ it  also holds for $k=j+1$. By setting $k=j$ in \eqref{eq_for_induction} we obtain that 

\begin{equation} \label{eq:inductive-Ek}
E_{j+1} \le E_j +  C' h^{s+1}(1+E_j)^{s+1}(L+M+1)^{s+1}[1 + h(1+E_j)(L+M+1)](E_j + E_{j+1} + 1).
\end{equation}
Using the assumption of induction in \eqref{ass_induc} we can obtain that $E_j\leq eE_0+1$ by setting $j=n$ in the right hand side. Using this inequality and the second condition in \eqref{eq:h-ineq}, we can write 
\begin{equation}
h\leq \frac{1}{C(eE_0+2)(L+M+1)} \leq \frac{1}{C(E_j+1)(L+M+1)}
\end{equation}
Using this expression we can simplify \eqref{eq:inductive-Ek} to
\begin{equation}\label{near_end_proof}
E_{j+1} \le E_j +  (1+{C}^{-1})C' h^{s+1}(1+E_j)^{s+1}(L+M+1)^{s+1} (E_j + E_{j+1} + 1).
\end{equation}
We can further show that 
\begin{align}\label{upp_bound_almost_done}
&(1+{C}^{-1})C' h^{s+1}(1+E_j)^{s+1}(L+M+1)^{s+1}\nonumber\\
& \leq (1+{C}^{-1})C' h^{s+1}(2+eE_0)^{s+1}(L+M+1)^{s+1} \leq \frac{1}{3N},
\end{align}
where the first inequality holds since $E_j\leq eE_0+1$ and the second inequality holds due to the second condition in \eqref{eq:h-ineq_2}. Simplifying the right hand side of \eqref{near_end_proof} using the upper bound \eqref{upp_bound_almost_done} leads to
\begin{equation}\label{eq:proof_444}
E_{j+1} \le E_j + \frac{1}{3N}(E_j + E_{j+1} + 1).
\end{equation}
Regroup the terms in \eqref{eq:proof_444} to obtain that $E_{j+1}$ is upper bounded by
\begin{equation}
E_{j+1} \le \left(\frac{1 + \frac{1}{3N}}{1-\frac{1}{3N}}\right)E_j + \frac{1}{3N-1}
\end{equation}
Now replace $E_j$ by its upper bound in \eqref{ass_induc} to obtain
\begin{align}\label{ttt_1}
E_{j+1} 
&\le \left(\frac{1 + \frac{1}{3N}}{1-\frac{1}{3N}}\right)\left(\left(1+\frac{1}{N}\right)^j E_0 + \frac{j}{N}\right) + \frac{1}{3N-1}\nonumber\\
&= \left(\frac{1 + \frac{1}{3N}}{1-\frac{1}{3N}}\right)\left(1+\frac{1}{N}\right)^j E_0 + \left(\frac{1 + \frac{1}{3N}}{1-\frac{1}{3N}}\right)\frac{j}{N} + \frac{1}{3N-1}\nonumber\\
&= \left(\frac{3N + 1}{3N-1}\right)\left(1+\frac{1}{N}\right)^j E_0 +
 \left(\frac{3N+1}{3N-1}\right)\frac{j}{N} + \frac{1}{3N-1}\nonumber\\
 &\le \left(1+\frac{1}{N}\right)^{j+1} E_0 +
 \left(\frac{3N+1}{3N-1}\right)\frac{j}{N} + \frac{1}{3N-1},
\end{align}
where the first inequality holds since $\frac{3N+1}{3N-1}\leq \frac{N+1}{N}$ and the last inequality follows from $1+\frac{2}{3N-1} \ge 1+\frac{1}{N}$. Further, we can show that
\begin{align}\label{ttt_2}
 \left(\frac{3N+1}{3N-1}\right)\frac{j}{N} + \frac{1}{3N-1} 
 &= \left(1+\frac{2}{3N-1}\right)\frac{j}{N} + \frac{1}{3N-1} \nonumber\\
  &= \frac{j}{N} + \left(\frac{2}{3N-1}\right)\frac{j}{N} + \frac{1}{3N-1} \nonumber\\
  &\leq  \frac{j}{N} + \left(\frac{2}{3N-1}\right)\frac{N-1}{N} + \frac{1}{3N-1} \nonumber\\
    &=  \frac{j}{N}+\frac{1}{N} \left(\frac{3N-2}{3N-1}\right)\nonumber\\
        &\leq  \frac{j+1}{N},
\end{align}
where in the first inequality we use the fact that $j\leq N-1$. Using the inequalities in \eqref{ttt_1} and \eqref{ttt_2} we can conclude that 
\begin{align}
E_{j+1}  \le \left(1+ \frac{1}{N}\right)^{j+1}E_0 + \frac{j+1}{N},
\end{align}

Therefore, the inequality in \eqref{claim_of_induction} is also true for $k=j+1$. The proof is complete by induction and we can write 
\begin{equation}
E_N \le eE_0 + 1.
\end{equation}
Now if we reconsider the conditions on $h$ in \eqref{eq:h-ineq} and \eqref{eq:h-ineq_2}, we can conclude that there exists a constant $\tilde{C}$ that is determined by $p, s$ and the numerical integrator, such that 
\begin{equation}
h \le \tilde{C} \frac{N^{-1/(s+1)}}{(L+M+1)(eE_0+1)},
\end{equation}
satisfies all the inequalities in  \eqref{eq:h-ineq} and \eqref{eq:h-ineq_2}.

\section{Bounding operator norms of derivatives and discretization errors of Lyapunov functions}

\begin{lemma} \label{lemma:bound-high-order-f}
Given state $y_c = [v_c; x_c; t_c]$ with $t_c \ge 1$ and the radius $R = \frac{1}{t_c}$, if $ B(x_c, R) \subseteq \mathcal{A}$ (defined in~\eqref{eq:set-A}) and hence Assumptions \ref{assump:growth},\ref{assump:differentiable} hold, then for all $ y \in U_{R, 0.2}(y_c)$ we can write
\begin{equation}\label{eq:claim_stay_in_the_ball}
\|\nabla^{(i)}f(x)\| \le p(M+L+1)\frac{\mathcal{E}(y_c) + 1}{t_c^{p-i}}.
\end{equation}
\end{lemma}

\begin{proof}
Based on Assumption \ref{assump:differentiable}, we know that 
\begin{align}\label{proof_lemma_10_eq_100}
\|\nabla^{(p)} f(x)\| \le M.
\end{align}
We further can show that the norm $\|\nabla^{(p-1)} f(x)\|$ is upper hounded by 
\begin{align}\label{proof_lemma_10_eq_200}
\|\nabla^{(p-1)} f(x)\| 
&= \|\nabla^{(p-1)} f(x_c) + \nabla^{(p-1)} f(x) - \nabla^{(p-1)} f(x_c)\| \nonumber\\
&\leq \|\nabla^{(p-1)} f(x_c)\| + \|\nabla^{(p-1)} f(x) - \nabla^{(p-1)} f(x_c)\|
\end{align}
Using the bound in \eqref{proof_lemma_10_eq_100} and the mean value theorem we can show that $\|\nabla^{(p-1)} f(x) - \nabla^{(p-1)} f(x_c)\|\leq M \|x-x_c\| \leq MR$, where the last inequality follows from $y \in U_{R, 0.2}(y_c)$. Applying this substitution into \eqref{proof_lemma_10_eq_200} implies that 
\begin{align}
\|\nabla^{(p-1)} f(x)\| 
&\leq  \|\nabla^{(p-1)} f(x_c)\| + MR \nonumber\\
& \le [L(f(x_c)-f(x^*))]^{\frac{1}{p}} + MR,
\end{align}
where the first inequality holds due to definition of operator norms and the last inequality holds due to the condition in Assumption \ref{assump:growth}. By following the same steps one can show that
\begin{align}
\|\nabla^{(p-2)} f(x)\| \le [L(f(x_c)-f(x^*))]^{\frac{2}{p}} + R[[L(f(x_c)-f(x^*))]^{\frac{1}{p}} + MR] 
\end{align}
By iteratively applying this procedure we obtain that if $y=[x;v;t]\in \reals^{2d+1}$ belongs to the set $U_{R, 0.2}(y_c)$, then we have
\begin{equation}\label{proof_eq_100}
\|\nabla^{(i)}f(x)\| \leq MR^{p-i} + \sum_{j=i}^{p-1} [L(f(x_c) - f(x^*))]^{\frac{p-j}{p}}R^{j-i}.
\end{equation}
Notice that since $\frac{p-j}{p}\leq 1$ for $j=1,\dots,p-1$, it follows that we can write $L^{\frac{p-j}{p}} \le 1+L$. Moreover, the definition of the Lyapunov function $\mathcal{E}$ in \eqref{eq:lyapunov} implies that 
\begin{equation}
[f(x_c) - f(x^*)]^{\frac{p-j}{p}} \le \frac{\mathcal{E}(y_c)^{\frac{p-j}{p}}}{t_c^{p-j}} \le \frac{1 + \mathcal{E}(y_c)}{t_c^{p-j}}
\end{equation}
where the last inequality follows from the fact that $\mathcal{E}(y_c)^{\frac{p-j}{p}} \le 1+\mathcal{E}(y_c)$ for $j=1,\dots,p-1$.
Therefore, we can simplify the upper bound in \eqref{proof_eq_100} by 
\begin{align} \label{eq:compute-Mi}
\|\nabla^{(i)}f(x)\| \leq MR^{p-i}+ \sum_{j=i}^{p} \frac{(1+L)(1+\mathcal{E}(y_c))}{t_c^{p-j}} R^{j-i}.
\end{align}
By replacing the radius $R$ with $1/t_c$ we obtain that 
\begin{align} \label{eq:compute-Mi_2}
\|\nabla^{(i)}f(x)\| 
& \leq \frac{M}{t_c^{p-i}}+ \sum_{j=i}^{p} \frac{(1+L)(1+\mathcal{E}(y_c))}{t_c^{p-i}} \nonumber\\
& = \frac{M+p(1+L)(1+\mathcal{E}(y_c))}{t_c^{p-i}}
\end{align}
As the Lyapunov function $\mathcal{E}(y_c)$ is always non-negative, we can write $M \leq M p(1+\mathcal{E}(y_c))$. Applying this substitution into \eqref{eq:compute-Mi_2} yields
\begin{align} \label{eq:compute-Mi_3}
\|\nabla^{(i)}f(x)\| 
 \leq  \frac{p(L+M+1)(1+\mathcal{E}(y_c))}{t_c^{p-i}},
\end{align}
and the claim in \eqref{eq:claim_stay_in_the_ball} follows.
\end{proof}

\begin{lemma} \label{lemma:bound-F}
If $ B(x_c, R) \subseteq \mathcal{A}$ (defined in~\eqref{eq:set-A}) and hence Assumptions \ref{assump:growth} and \ref{assump:differentiable} hold, there exists a constant $C$ determined by $p$ such that, $\forall y \in U_{R, 0.2}(y_c)$ where $y_c = [v_c; x_c; t_c]$, $t_c \ge 1$ and $R = \frac{1}{t_c}$,  we have
\begin{align}
\|\pi_{x,v} F(y)\|= \le \frac{C(\mathcal{E}(y_c) + 1)(L+M+1)}{t_c}.
\end{align}
\end{lemma}
\begin{proof}
According to Lemma \ref{lemma:bound-high-order-f}, we can write that 
\begin{align}\label{eq:proof_int_lemma_50}
\|\nabla f(x)\| \le p(M+L+1)\frac{\mathcal{E}(y_c) + 1}{t_c^{p-1}}.
\end{align}
Further, the definition of the Lyapunov function in \eqref{eq:lyapunov} implies that
\begin{equation}\label{eq:proof_int_lemma_100}
\|v_c\| \le \frac{2p\mathcal{E}(y_c)^{0.5}}{t_c}.
\end{equation}
Since $y \in U_{R, 0.2}(y_c)$, we have that 
\begin{align}\label{eq:proof_int_lemma_200}
|t-t_c|\le 0.2,\qquad  \|v-v_c\|\le R, \qquad \|x-x_c\|\le R.
\end{align}
Further, based on the dynamical system in \eqref{eq:dynamics}, we can write
\begin{align}\label{eq:proof_int_lemma_300}
\|\pi_{x,v} F(y)\|
&= 
\left\|\begin{bmatrix}
-\frac{2p+1}{t} v - p^2t^{p-2} \nabla f(x)\\
v\\
\end{bmatrix}\right\| \nonumber\\
& \le \frac{2p+1}{t}\| v\| +\|p^2t^{p-2} \nabla f(x)\| +\|v\|\nonumber\\
& \le \left(\frac{2p+1}{t}+1\right)(\| v_c\| +\|v_c-v\|) +p^2t^{p-2} \|\nabla f(x)\| ,
\end{align}
where the first inequality is obtained by using the property of norm, and in the last one we use the triangle inequality. Note that according to \eqref{eq:proof_int_lemma_200} we have $t\geq t_c -0.2 $. Since $t_c\geq 1$ it implies that  $t\geq 0.8 t_c$. In addition we can also show that $t\leq t_c+0.2\leq 1.2 t_c$. Applying these bounds into \eqref{eq:proof_int_lemma_300} yields
\begin{align}\label{eq:proof_int_lemma_400}
\|\pi_{x,v} F(y)\|
\le \left(\frac{p+1}{0.8 t_c}+1\right)(\| v_c\| +\|v_c-v\|) +(1.2)^{p-2}p^2t_c^{p-2} \|\nabla f(x)\|
\end{align}
Replace $\|\nabla f(x)\|$, $\| v_c\|$, and $\|v_c-v\|$ in \eqref{eq:proof_int_lemma_400} by their upper bounds in \eqref{eq:proof_int_lemma_50}, \eqref{eq:proof_int_lemma_100}, and \eqref{eq:proof_int_lemma_200}, respectively, to obtain
\begin{align}\label{eq:proof_int_lemma_500}
\|\pi_{x,v} F(y)\|
&\le \left(\frac{p+1}{0.8 t_c}+1\right)\left( \frac{2p\mathcal{E}(y_c)^{0.5}}{t_c} +R\right) +(1.2)^{p-2}p^3 (M+L+1)\frac{\mathcal{E}(y_c) + 1}{t_c}\nonumber\\
&\le \left(\frac{p+1}{0.8 t_c}+1\right)\left( \frac{2p(\mathcal{E}(y_c)+1)+1}{t_c} \right) +(1.2)^{p-2}p^3 (M+L+1)\frac{\mathcal{E}(y_c) + 1}{t_c},
\end{align}
where in the second inequality we replace $R$ by $1/t_c$ and $\mathcal{E}(y_c)^{0.5}$ by its upper bound $\mathcal{E}(y_c)+1$. Considering that $t_c\geq 1$ and the result in \eqref{eq:proof_int_lemma_500} w obtain that there exists a constant $C$ such that 
\begin{align}
\|\pi_{x,v} F(y)\|
\le \frac{C(\mathcal{E}(y_c) + 1)(L+M+1)}{t_c},
\end{align}
where $C$ only depends on $p$.
\end{proof}

\begin{lemma} \label{lemma:bound-high-order-F}
Given state $y_c = [v_c, x_c, t_c]$ with $t_c \ge 1$, let $R = \frac{1}{t_c}$. If $ B(x_c, R) \subseteq \mathcal{A}$ (defined in~\eqref{eq:set-A}) and hence Assumptions \ref{assump:growth},\ref{assump:differentiable} hold, then when $h \le  \min\{0.2,\frac{1}{(1+\kappa)C(\mathcal{E}(y_c) + 1)(L+M+1)}\}$, we have
\begin{align}
\left\| \frac{\partial^q \varphi_h(y_c)}{\partial h^q} \right\| \le  \frac{C_0 [\mathcal{E}(y_c) + 1]^{q}(L+M+1)^q}{t_c}, 
\end{align}
and 
\begin{align}
\left\| \frac{\partial^q \Phi_h(y_c)}{\partial h^q} \right\| \le  \frac{C_1 [1 + \mathcal{E}(y_c)]^q (L+M+1)^q + C_2 h[1+\mathcal{E}(y_c)]^{q+1}(L+M+1)^{p+1}}{t_c},
\end{align}
where $C$ and $\kappa$ are the same constants as in Lemma \ref{lemma:bound-F}. Further, the constants $C_1,C_2, C_3$ are determined by $p$, $q$, and the integrator.
\end{lemma}


\begin{remark}
In the proof below, we reuse variants of symbol $C$(e.g.$C_1, C_2, \tilde{C}$) to hide constants determined by $p, q$ and the integrator. We recommend readers to focus on the degree of the polynomials in $(L+M+1), \mathcal{E}(y_c), h, t_c$, and check that the rest can be upper-bounded by variants of symbol $C$.
We frequently use two tricks in this section. First, for $a \in (0, 1)$, we can bound
\begin{equation}
c^a \le c+1
\end{equation}
Second, note that given $t_c \ge 1$,for any $n>0$, there exist constants $C_1, C_2,C_3$ determined by $n$ such that for all $t$ subject to $|t-t_c| \le 0.2$, 
\begin{equation}
\frac{1}{t^n} \le \frac{C_1}{t_c^n} \le C_2\\
t^n \le C_3 t_c^n
\end{equation}
\end{remark}

\begin{proof}	
Notice that the system dynamic function $F:\reals^{2d+1} \to \reals^{2d+1}$ in Equation (\ref{eq:dynamics}) is a vector valued multivariate function. We denote its $i_{th}$ order derivatives by $\nabla^{(i)}F(y)$, which is a $\repeatn{$(2d+1) \times$}{$\times (2d+1)$}{i+1\ times}$ tensor. The tensor is symmetric by continuity and Schwartz theorem. As a shorthand, we use $\nabla^{(i)}F$ to denote $\nabla^{(i)}F(y)$. We know that
$y^{(i)} = F^{(i-1)}(y) = \frac{\partial^i y}{\partial t^i}$. Notice that $F^{(i-1)}(y)$ is a vector. As an example, we can write
\begin{align}
y^{(1)} =& F\nonumber\\
y^{(2)} =& F^{(1)} = \nabla F(F)\nonumber \\
y^{(3)} =& F^{(2)} = \nabla^{(2)}F(F, F) + \nabla F(\nabla F(F)) .
\end{align}
The derivative $\nabla^{(i)}F(y)$ can be interpreted as a linear map: $\nabla^{(i)}F: \repeatn{$\reals^{2d+1} \times$}{$\times \reals^{2d+1}$}{i \ \text{times}} \to \reals^{2d+1}$. $\nabla^{(2)}F(F_1, F_2)$ maps $F_1, F_2$ to some element in $\reals^{2d+1}$.  Enumerating the expressions will soon get very complicated. However, we can express them compactly with elementary differentials summarized in Appendix \ref{sec:elementary} (see Chapter 3.1 in \cite{hairer-textbook} for details). 

First we bound $\nabla^{(i)}F$ by explicitly computing its entries. Let $a(t) = p^2t^{p-2}$ and $b(t) = \frac{2p+1}{t}$. Based on the definition in \eqref{eq:dynamics}, we obtain that
\begin{align}
 \frac{\partial^{k + 1} F}{ \partial v \partial t^k } =& \begin{bmatrix} -b^{(k)}(t) I \\ I^{(k)}  \\ 0\\ \end{bmatrix},  \qquad
\frac{\partial^{k} F}{\partial t^k  } = \begin{bmatrix} -b^{(k)}(t) v - a^{(k)}(t) \nabla f(x) \\ 0\\ 0\\ \end{bmatrix}, \nonumber\\
\frac{\partial^{i + k} F}{\partial x^i \partial t^k} =& \begin{bmatrix}
-a^{(k)}(t) \nabla^{i+1} f(x)\\0\\0 \end{bmatrix}, \qquad
\frac{\partial^{i} F}{\partial x^i  } = \begin{bmatrix} -a(t) \nabla^{i+1} f(x) \\ 0\\ 0\\ \end{bmatrix},\nonumber\\
 \frac{\partial^{i + j} F}{ \partial v^j \partial x^i} =& 0, \qquad  \frac{\partial F}{ \partial v } =\begin{bmatrix} \frac{2p+1}{t}I \\ I\\ 0\\ \end{bmatrix},
 \qquad \frac{\partial^{ j } F}{ \partial v^j }= 0, j \ge 2. \\
\end{align}
 For any vector $y = [v; x; t] \in U_{R, 0.2}(y_c)$, we can show that the norm of $\nabla^{(n)}F$ is upper bounded by
\begin{align} \label{eq:bound-grad-F-components}
\|\nabla^{(n)}F(F_1, F_2, ..., F_n)\| & \le \|a(t) \nabla^{(n+1)}f(x)\|\prod_{i \in [n]} \|\pi_x F_i\| \nonumber\\
& + \|b^{(n)}(t)v + a^{(n)}(t)\nabla f(x)\|\prod_{i \in [n]} \|\pi_t F_i\| \nonumber\\
& + \sum_{k \ge 1}^{n-1} \sum_{\substack{S \subset [n]\\ |S|=k}}\|a^{(k)}(t)\nabla^{(n-k+1)}f(x)\|\left[\prod_{s \in S } \|\pi_t F_s\| \right] \left[\prod_{s' \in [n]/S} \|\pi_x F_{s'}\|\right] \nonumber\\
& + \sum_{i \in [n]}\|b^{(n-1)}(t)+1\|\|\pi_v F_i\|\prod_{j \neq i}\|\pi_t F_j\|.
\end{align}

Using the definition of the Lyapunov function $\mathcal{E}$ and the definition of the set $U_{R, 0.2}(y_c)$ it can be shown that
\begin{align}\label{proof_last_lemma_100}
\|v_c\| \le \frac{\mathcal{E}(y_c)^{0.5}}{t_c} \le \frac{\mathcal{E}(y_c) + 1}{t_c},\quad t_c \ge 1,\quad |t-t_c| \le 0.2, \quad \|v-v_c\| \le R.
\end{align} 
Further, the result in Lemma \ref{lemma:bound-high-order-f} implies that
\begin{equation}\label{proof_last_lemma_200}
\|\nabla^{(i)}f(x)\| \le p(M+L+1)\frac{\mathcal{E}(y_c) + 1}{t_c^{p-i}}.
\end{equation}
Substituting the upper bounds in \eqref{proof_last_lemma_100} and \eqref{proof_last_lemma_200} into \eqref{eq:bound-grad-F-components} implies that for $n=1,\dots,p$ we can write
\begin{align} \label{eq:bound-grad-Fi}
&\|\nabla^{(n)}F(F_1, F_2, ..., F_n)\| \nonumber\\
& \le C_1(M+L+1)[\mathcal{E}(y_c) + 1]t_c^{n-1} \prod_{i \in [n]}\|\pi_x F_i\|  \nonumber\\
&\quad+ C_2 (M+L+1)\left[\mathcal{E}(y_c) + 1\right]t_c^{-n-1}\prod_{i \in [n]}\|\pi_t F_i\| \nonumber\\
&\quad+ C_3 (M+L+1) \sum_{k \ge 1}^{p-1} \left[\mathcal{E}(F_c) + 1\right] t_c^{n-2k-1} \sum_{\substack{S \subset [n]\\ |S|=k}}\left[\prod_{s \in S } \|\pi_t F_s\|\right]\left[\prod_{s' \in [n]/S} \|\pi_x F_{s'}\|\right]  \nonumber\\
& \quad+ C_4 \sum_{i \in [n]}\left[1 + \frac{1}{t_c^n}\right] \|\pi_v F_i\|\prod_{j \neq i}\|\pi_t F_j\|,
\end{align}
where $C_1, C_2, C_3$, and $C_4$ only depend on $n$ and $p$. 

For $n = p, p+1, ..., s$, we can get similar bounds. To do so, not only we use the result in \eqref{proof_last_lemma_200}, but also we use the bounds guaranteed by Assumption \ref{assump:differentiable}. Hence, for $n = p, p+1, ..., s$ it holds
\begin{align} \label{eq:bound-grad-Fi-p}
&\|\nabla^{(n)}F(F_1, F_2, ..., F_n)\| \nonumber\\
& \le C_1 Mt_c^{p-2} \prod_{i \in [n]}\|\pi_x F_i\|  \nonumber\\
&\quad+ C_2 (M+L+1)[\mathcal{E}(y_c) + 1]t_c^{-n-1}\prod_{i \in [n]}\|\pi_t F_i\| \nonumber\\
&\quad+ C_3 \sum_{k \ge 1}^{p-1} (M+L+1)[\mathcal{E}(y_c) + 1]t_c^{p-k-2} \sum_{\substack{S \subset [n]\\ |S|=k}}\left[\prod_{s \in S } \|\pi_t F_s\|\right]\left[\prod_{s' \in [n]/S} \|\pi_x F_{s'}\|\right]  \nonumber\\
&\quad + C_4 \sum_{i \in [n]} \left[1 + \frac{1}{t_c^n}\right] \|\pi_v F_i\|\prod_{j \neq i}\|\pi_t F_j\|.
\end{align}

Finally we are ready to bound the time derivatives. We first bound the elementary differentials $F(\tau)$ defined in Section \ref{sec:elementary} Definition \ref{def:elementary}. Let $F(\tau) = F(\tau)(y)$ for convenience. We claim that when $|\tau| \le q$, then $\forall y \in U_{R, 0.2}(y_c)$
\begin{align}  \label{eq:induct-assump}
\|\pi_t F(\tau)\| \le 1,\qquad \|\pi_{v,x} F(\tau) \| \le C_{|\tau|} (L+M+1)^{|\tau|} \frac{[\mathcal{E}(y_c) + 1]^{|\tau|}}{t_c},
\end{align}
where the constant $C_q$ only depends on $p$ and $q$. 
We use induction to prove the claims in \eqref{eq:induct-assump}. The base case is trivial as we have shown  in Lemma \ref{lemma:bound-F} that $\|\pi_{x,v}F(\bullet)(y)\| = \|\pi_{x,v} F(y)\| \le \frac{C(\mathcal{E}(y_c)+1)(L+M)}{t_c}$, and $\|\pi_{t}F(\bullet)(y)\|=\|\pi_t F(y)\| = 1$. Since the last coordinate grows linearly with rate 1 no matter what $x, v$ are, it can be shown that
\begin{align}
\pi_t F(\tau)(y) = 0, \forall |\tau| \ge 2.
\end{align}
We hence focus on proving the upper bound for the norm $\|\pi_{x,v}F(\tau)(y)\|$ in \eqref{eq:induct-assump}.

Now assume $|\tau| = q$ and it has $m$ subtrees attached to the root, $\tau = [\tau_1, ..., \tau_m]$ with $\sum_{i=1}^m |\tau_i| = q - 1$. When $m \le p-1$, by \eqref{eq:bound-grad-Fi} we obtain
\begin{align} \label{eq:F-induc}
&\|\nabla^{(m)}F(F(\tau_1),  ..., F(\tau_m))\| \nonumber\\
& \le C_1 [(M+L+1) (\mathcal{E}(y_c) + 1)]t_c^{m-1} \prod_{i \in [m]}\|\pi_x F(\tau_i)\| \nonumber\\
&\quad + C_2 (M+L+1)[\mathcal{E}(y_c) + 1]t_c^{-m-1}\prod_{i \in [m]}\|\pi_t F(\tau_i)\|\nonumber\\
&\quad + C_3 \sum_{k \ge 1}^{m-1} [(M+L+1) (\mathcal{E}(y_c) + 1)1] t_c^{m-2k-1} \sum_{\substack{S \subset [m]\\ |S|=k}}\left[\prod_{s \in S } \|\pi_t F(\tau_s)\|\right]\left[\prod_{s' \in [m]/S} \|\pi_x F(\tau_{s'})\|\right] \nonumber\\
& \quad + C_4 \sum_{i \in [m]} \left[1 + \frac{1}{t_c^n}\right] \|\pi_v F(\tau_i)\|\prod_{j \neq i}\|\pi_t F(\tau_j)\|.
\end{align}
Notice that $|\tau_i| \le q-1$. By inductive assumption in \eqref{eq:induct-assump} we can write
\begin{align} 
\|\pi_t F(\tau_i)\| &\le 1 \quad \text{for all} \ i=1\dots,m\\
\prod_{i \in S}\|\pi_{v,x} F(\tau_i)\| & \le C_n (L+M+1)^n \frac{[\mathcal{E}(y_c) + 1]^n}{t_c^{|S|}}, \quad \text{where} \ n = \sum_i |\tau_i|.
\end{align}
Apply these substitutions into \eqref{eq:F-induc} to and use the inequality $\sum_i |\tau_i|\leq q-1$ to obtain that
\begin{align} \label{eq:F-induc_22}
\|\nabla^{(m)}F(F(\tau_1),  ..., F(\tau_m))\|  \le C_q \frac{[\mathcal{E}(y_c) + 1]^q(M+L+1)^q}{t_c}.
\end{align}
Hence, since $\|\pi_{x,v} F(\tau) \| \le \|\nabla^{(m)}F(F(\tau_1),  ..., F(\tau_m))\|$ we obtain that 
\begin{align} 
\|\pi_{x,v} F(\tau) \|  \le C_q \frac{[\mathcal{E}(y_c) + 1]^q(M+L+1)^q}{t_c}.
\end{align}
Similarly, for $m \ge p$, by (\ref{eq:bound-grad-Fi-p}) we can write
\begin{align}\label{kashk}
& \|\nabla^{(m)}F(F(\tau_1),  ..., F(\tau_m))\| \nonumber\\
&\le C_1 Mt_c^{p-2} \prod_{i \in [m]}\|\pi_x F(\tau_i)\| \nonumber\\
&\quad + C_2 (M+L+1)[\mathcal{E}(y_c) + 1]t^{-n-1}\prod_{i \in [n]}\|\pi_t F(\tau_i)\|\nonumber\\
&\quad + C_3 \sum_{k \ge 1}^{m-1} [(M+L+1) (\mathcal{E}(y_c) + 1)1] t_c^{p-k-2} \sum_{\substack{S \subset [m]\\ |S|=k}}\left[\prod_{s \in S } \|\pi_t F(\tau_s)\|\right]\left[\prod_{s' \in [m]/S} \|\pi_x F(\tau_{s'})\|\right] \nonumber\\
& \quad + C_4 \sum_{i \in [m]}\left[1 + \frac{1}{t_c^n}\right] \|\pi_v F(\tau_i)\|\prod_{j \neq i}\|\pi_t F(\tau_j)\|.
\end{align}
Plug in the induction assumption in $\eqref{eq:induct-assump}$ into \eqref{kashk} to obtain
\begin{align} \label{eq:bound-elementary}
\|\pi_{x,v} F(\tau) \| \le \|\nabla^{(m)}F(F(\tau_1),  ..., F(\tau_m))\| \le C_q \frac{[\mathcal{E}(y_c) + 1]^q(M+L+1)^q}{t_c}.
\end{align}
Hence, the proof is complete by induction. 

Now we  proceed to derive an upper bound for higher order time derivatives. By Lemma \ref{lemma:high-derivative-solution} we can write
$$\| \frac{\partial^q \varphi_h(y_c)}{\partial h^q} \| =   \|F^{(q-1)}(\varphi_h(y_c))\| = \|\sum_{|\tau|=q}\alpha(\tau)F(\tau)(\varphi_h(y_c))\|.$$

By Lemma \ref{lemma:bound-F}, we know that when $h \le  \min\{0.2,\frac{1}{(1+\kappa)C(\mathcal{E}(y_c) + 1)(M+L)}\}$, $y \in U_{R, 0.2}(y_c)$. Therefore, (\ref{eq:bound-elementary}) holds. Hence, there exists a constant $C$ determined by $p,q$ such that

$$\| \frac{\partial^q \varphi_h(y_c)}{\partial h^q} \| \le  \frac{C [\mathcal{E}(y_c) + 1]^{q}(M+L+1)^q}{t_c}. $$

Similarly by Lemma \ref{lemma:high-derivative-numerical}, we have the following equation
$$\frac{\partial^q \Phi_h(y_c)}{\partial h^q} = \sum_{i \le S} b_i [h \frac{\partial^q F(g_i)}{\partial h^q} + q\frac{\partial^{q-1} F(g_i)}{\partial h^q}]$$
Here, $\frac{\partial^q F(g_i)}{\partial h^q}$ has the same recursive tree structure as $F^{(q)}(y)$, except that we need to replace all $F$ in the expression by 
$\frac{\partial g_i}{\partial h}$ and all $\nabla^{(n)}F(y)$ by $\nabla^{(n)}F(g_i)$. By Definition \ref{def:rk} and Lemma \ref{lemma:bound-F}, we know that 
$$\|\frac{\pi_{x,v}\partial g_i}{\partial h}\| \le \sum_{j\le i-1} |a_{ij}| \frac{C(\mathcal{E}(y_c)+1)(M+L+1)}{t_c} ,\qquad \|\frac{\pi_{t}\partial g_i}{\partial h}\| = |\sum_{j\le i-1} a_{ij}|. $$
We also know by lemma \ref{lemma:bound-F} that $\forall i, g_i \in U_{R, 0.2}(y_c)$. Hence the bounds for $\|\nabla^{(n)}F(y)\|$ also holds for $\nabla^{(n)}F(g_i)$. Therefore, by the same argument as for bounding $\| \frac{\partial^q \varphi_h(y_c)}{\partial h^q} \|$, we will get same bounds for  $\|\frac{\partial^q F(g_i)}{\partial h^q}\|$ up to a constant factor determined by the integrator. Based on this, we conclude  that
$$\| \frac{\partial^q \Phi_h(y_c)}{\partial h^q} \| \le  \frac{C [(L+M+1)(1 + \mathcal{E}(y_c))]^q + C' h[(L+M+1)(1 + \mathcal{E}(y_c))]^{(q+1)}}{t_c}, $$
where the constants are determined by $p,q$ and the integrator.
\end{proof}

\begin{lemma} \label{lemma:bound-discrete-lyapunov}
Suppose the conditions in Proposition ~\ref{prop:discrete-error} hold. Then, we have that
\begin{align}\label{last_lemma_claim}
&\|\mathcal{E}(\Phi_h(y_k)) - \mathcal{E}(\varphi_h(y_k))\| \nonumber\\
&\le C h^{s+1}[(1+E_k)^{s+1}(L+M+1)^{s+1} + h(1+E_k)^{s+2}(L+M+1)^{s+2}](E_k + E_{k+1} + 1),
\end{align}
where $C$ only depends on $p, s$ and the numerical integrator. 
\end{lemma}
\begin{proof}
Denote $\hat{y} = \Phi_h(y_k), \tilde{y} = \varphi_h(y_k)$. Notice that $\tilde{t} = \hat{t} = t_k + h$.  In fact, because we start the simulation at $t_c=1$ and we require that $h \le 0.2$, we have  
\begin{align} \label{eq:bound-tk-t}
\frac{t_k}{\tilde{t}} = \frac{t_k}{t_k + h} \in \left[\frac{5}{6} \ ,\ 1\right].
\end{align}
Now using the definition of the Lyapunov function $\mathcal{E}$ we can show that  
\begin{align}\label{eq:diff-obj}
\|\mathcal{E}(\hat{y}) - \mathcal{E}(\tilde{y})\| \le& \frac{\tilde{t}^2}{4p^2}\left|\|\tilde{v}\|^2 - \|\hat{v}\|^2\right| + \left|\|\tilde{x} + \frac{\tilde{t}}{2p}\tilde{v} - x^*\|^2 - \| \hat{x} + \frac{\hat{t}}{2p}\hat{v} - x^* \|^2\right| + \tilde{t}^p(|f(\tilde{x}) - f(\hat{x})|) \nonumber\\
\le& \frac{2\hat{t}^2}{4p^2}(\|\tilde{v}-\hat{v}\|\|\tilde{v}+\hat{v}\|) + \tilde{t}^p (\|\tilde{x} - \hat{x}\|)(\|\nabla f(\tilde{x})\|+\|\nabla f(\hat{x})\|)\nonumber\\
&+ 2\left\|\tilde{x} - \hat{x} + \frac{\tilde{t}}{2p}(\tilde{v} - \hat{v})\right\| \left\|\tilde{x} + \frac{\tilde{t}}{2p}\tilde{v} - x^* +\hat{x} + \frac{\hat{t}}{2p}\hat{v} - x^*\right\|,
\end{align}

where to derive the second inequality we used the convexity of the function $f$ which implies

\begin{align}
\inner{y-x, \nabla f(y)} \le f(x) - f(y) \le \inner{x-y, \nabla f(x)}.
\end{align}
Recall that $E_k = \mathcal{E}(y_k)$, $E_{k+1} = \mathcal{E}(\hat{y}) = \mathcal{E}(\Phi_h(y_k))$, $\tilde{E}_{k+1} = \mathcal{E}(\tilde{y}) = \mathcal{E}(\varphi_h(y_k))$. According to Proposition~\ref{prop_dec_energy} we know that $\tilde{E}_{k+1} \le E_k$, and therefore $\tilde{E}_{k+1}$ is upper bounded by $E_k$. Therefore, we can write
\begin{align}\label{first_bounds_ineq}
&\|\tilde{v}\| \le \frac{\sqrt[]{\tilde{E}_{k+1}}}{\tilde{t}} \le \frac{\sqrt[]{E_k}}{\tilde{t}} \le \frac{E_k+1}{\tilde{t}},\qquad \|\hat{v}\| \le \frac{E_{k+1}+1}{\hat{t}},\nonumber\\
&\left\|\tilde{x} + \frac{\tilde{t}}{2p}\tilde{v} - x^*\right\| \le \sqrt[]{E_k} \le E_k + 1,\qquad
\left\|\hat{x} + \frac{\hat{t}}{2p}\hat{v} - x^*\right\| \le E_{k+1} + 1.
\end{align}
Further, by Assumption \ref{assump:growth}, we have that
\begin{align}\label{second_bounds_ineq}
\|\nabla f (\tilde{x})\| \le \frac{L(E_{k}+1)}{\tilde{t}^{p-1}}, \qquad \|\nabla f (\hat{x})\| \le L(f(\hat{x})-f(x^*))^{\frac{p-1}{p}} \le L(\frac{E_{k+1}}{\hat{t}^p})^{\frac{p-1}{p}} \le \frac{L(E_{k+1}+1)}{\hat{t}^{p-1}}.
\end{align}
In addition, by Proposition~\ref{prop:discrete-error}, we know that for some constant C determined by $p, s, L, M$ and the integrator, it holds
\begin{align}\label{third_bounds_ineq}
&\max\{\|\tilde{v} - \hat{v}\|, \|\tilde{x} - \hat{x}\|\}\nonumber\\
& \le C h^{s+1}\left[ \frac{[1+\mathcal{E}(y_k)]^{s+1}(L+M+1)^{s+1}}{t_k}\ +\ h\ \frac{[1+\mathcal{E}(y_k)]^{s+2}(L+M+1)^{s+2}}{t_k}\right].
\end{align}
Define $\mathcal{M}:=[ \frac{[1+\mathcal{E}(y_k)]^{s+1}(L+M+1)^{s+1}}{t_k} + h \frac{[1+\mathcal{E}(y_k)]^{s+2}(L+M+1)^{s+2}}{t_k}]$. Use the upper bounds in \eqref{first_bounds_ineq}-\eqref{third_bounds_ineq} and the definition of $\mathcal{M}$ to simplify the right hand side of \eqref{eq:diff-obj} to 
\begin{align}
\|\mathcal{E}(\hat{y}) - \mathcal{E}(\tilde{y})\| \le & \frac{2\tilde{t}^2}{4p^2}Ch^{s+1}\mathcal{M}\frac{E_k+E_{k+1}+2}{\tilde{t}} 
+ \tilde{t}^p Ch^{s+1}\mathcal{M}\frac{L(E_{k+1}+E_k+2)}{\tilde{t}^{p-1}}
\nonumber \\
&\quad 
+ 2\left(1+\frac{t_k}{2p}\right)Ch^{s+1}\mathcal{M} (E_k + E_{k+1} + 2).
\end{align}
Now use the fact that $\frac{t_k}{\tilde{t}}$ is bounded by a constant as shown \eqref{eq:bound-tk-t}. Further, upper bound all the constants determined by $s, p$ and the numerical integrator, we obtain that
\begin{align}
&\|\mathcal{E}(\hat{y}) - \mathcal{E}(\tilde{y})\| \nonumber\\
&\le C' h^{s+1}[(1+E_k)^{s+1}(L+M+1)^{s+1} + h(1+E_k)^{s+2}(L+M+1)^{s+2}](E_k + E_{k+1} + 1),
\end{align}
and the claim in \eqref{last_lemma_claim} follows. 
\end{proof}

\section{Elementary differentials} \label{sec:elementary}
\begin{figure}[t]
    \centering
    \includegraphics[width=0.5\textwidth]{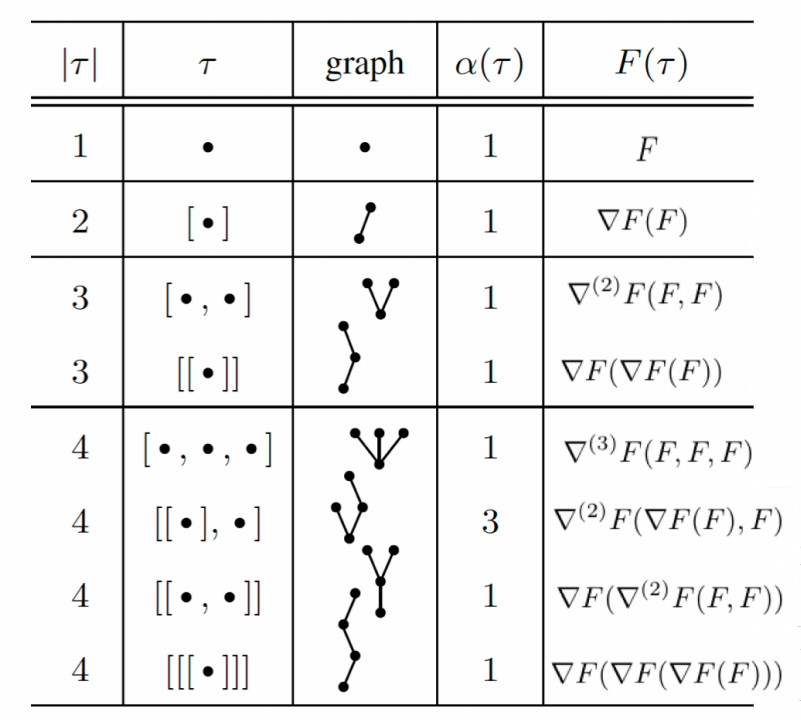}
    \caption{A figure adapted from \cite{hairer-textbook}. Example tree structures and corresponding function derivatives.}
    \label{fig:bseries}
\end{figure}

We briefly summarize some key results on elementary differentials from \cite{hairer-textbook}. For more details, please refer to chapter 3 of the book. Given a dynamical system 
$$\dot{y} = F(y)$$
we want to find a convenient way to express and compute its higher order derivatives. To do this, let $\tau$ denote a tree structure as illustrated in Figure \ref{fig:bseries}. $|\tau|$ is the number of nodes in $\tau$. Then we can adopt the following notations as in \cite{hairer-textbook}

\begin{definition} \label{def:elementary}
For a tree $\tau$, the elementary differential is a mapping $F(\tau): \reals^d \to \reals^d$, defined recursively by $F(\bullet)(y) = F(y)$ and 
$$F(\tau) = \nabla^{(m)}F(y)(F(\tau_1)(y), ..., F(\tau_m)(y))$$
for $\tau = [\tau_1, ..., \tau_m]$. Notice that $\sum_{i=1}^m |\tau_i| = |\tau| - 1$.
\end{definition}

Some examples are shown in Figure \ref{fig:bseries}. With this notation, the following results from \cite{hairer-textbook} Chapter 3.1 hold. The proof follows by recursively applying the product rule.
\begin{lemma} \label{lemma:high-derivative-solution}
The qth order derivative of the exact solution to $\dot{y} = F(y)$ is given by
$$y^{(q)}(t_c) = F^{(q-1)}(y_c) = \sum_{|\tau|=q}\alpha(\tau)F(\tau)(y_c)$$
for $y(t_c) = y_c$. $\alpha(\tau)$ is a positive integer determined by $\tau$ and counts the number of occurrences of the tree pattern $\tau$.
\end{lemma}
The next result is obtained by Leibniz rule. The expression for $\frac{\partial^q F(g_i)}{\partial h^q}$ can be calculated the same way as in Lemma \ref{lemma:high-derivative-solution}.
\begin{lemma} \label{lemma:high-derivative-numerical}
For a Runge-Kutta method defined in definition \ref{def:rk}, if $F$ is $q_{th}$ differentiable, then
\begin{align}\label{eq:Leibniz}
\frac{\partial^q \Phi_h(y_c)}{\partial h^q} = \sum_{i \le S} b_i [h \frac{\partial^q F(g_i)}{\partial h^q} + q\frac{\partial^{q-1} F(g_i)}{\partial h^q}]
\end{align}
where $\frac{\partial^q F(g_i)}{\partial h^q}$ has the same structure as $F^{(q)}(y)$ in lemma \ref{lemma:high-derivative-solution}, except that we need to replace all $F$ in the expression by 
$\frac{\partial g_i}{\partial h}$ and all $\nabla^{(n)}F(y)$ by $\nabla^{(n)}F(g_i)$.
\end{lemma}

\end{document}